\input ./style/arxiv-general.cfg
\documentclass[aap,MSNbibl,seceqn,dvips]{arximspdf}
\makeatletter
   \@ifpackageloaded{graphicx}{}{\usepackage{graphicx}}
\makeatother

\doi{10.1214/14-AAP1046} 
\volume{25}
\issue{4}
\pubyear{2015}
\firstpage{2168}
\lastpage{2214}
\docsubty{FLA}

\makeatletter
\newcommand{\lleft}{\left}
\newcommand{\rrvert}{\vert}
\newcommand{\rright}{\right}
\newcommand{\rrVert}{\Vert}
\newcommand{\llvert}{\vert}
\newcommand{\llVert}{\Vert}
\newtheorem{teo}{Theorem}[section]
\newtheorem{lem}[teo]{Lemma}
\newtheorem{cor}[teo]{Corollary}
\newproclaim{rem}[teo]{Remark}
\newproclaim{defn}[teo]{Definition}
\newproclaim{assum}[teo]{Assumption}
\def\by{\mathbf{y}}
\def\cF{\mathcal{F}}
\def\cG{\mathcal{G}}
\def\cX{\mathcal{X}}
\def\cY{\mathcal{Y}}
\def\cZ{\mathcal{Z}}
\def\essinf{\mathop{\operatorname{essinf}}}
\def\esssup{\mathop{\operatorname{esssup}}}
\def\liminf{\mathop{\underline{\operatorname{lim}}}}
\makeatother

\begin{document}
\begin{frontmatter}

\title{On well-posedness of forward--backward\\ SDEs---A unified approach}
\runtitle{Well-posedness of FBSDEs}

\begin{aug}
\author[A]{\fnms{Jin}~\snm{Ma}\corref{}\ead[label=e1]{jinma@usc.edu}\thanksref{T1}},
\author[B]{\fnms{Zhen}~\snm{Wu}\ead[label=e2]{wuzhen@sdu.edu.cn}\thanksref{T2}},
\author[B]{\fnms{Detao}~\snm{Zhang}\ead[label=e3]{zhangdetao@sdu.edu.cn}\thanksref{T3}}
\and
\author[A]{\fnms{Jianfeng}~\snm{Zhang}\ead[label=e4]{jianfenz@usc.edu}\thanksref{T4}}
\runauthor{Ma, Wu, Zhang and Zhang}
\affiliation{University of Southern California, Shandong University,\\
Shandong University and University of Southern California}
\address[A]{J.~Ma\\
J. Zhang\\
Department of Mathematics\\
University of Southern California\\
Los Angeles, California 90089\\
USA\\
\printead{e1}\\
\phantom{E-mail: }\printead*{e4}}
\address[B]{Z. Wu\\
D. Zhang\\
School of Mathematics\\
Shandong University\\
Jinan, 250100\\
P.R. China\\
\printead{e2}\\
\phantom{E-mail: }\printead*{e3}}
\end{aug}
%
\thankstext{T1}{Supported in part by NSF Grants DMS-0806017 and DMS-1106853.}
\thankstext{T2}{Supported in part by Chinese NSF Grants \#10921101, \#61174092 and the National Science Fund for Distinguished Young Scholars
of China \#11125102.}
\thankstext{T3}{Supported in part by Independent Innovation Foundation
of Shandong University (2011GN018).}
\thankstext{T4}{Supported in part by NSF Grant DMS-1008873.}

\received{\smonth{8} \syear{2013}}
\revised{\smonth{4} \syear{2014}}

%
\begin{abstract}
In this paper, we study the well-posedness of the Forward--Backward
Stochastic Differential Equations (FBSDE) in a general non-Markovian
framework. The main purpose is to find a unified scheme which combines
all \mbox{existing} methodology in the literature, and to address
some fundamental longstanding problems for non-Markovian
FBSDEs. An important device is a
\textit{decoupling random field} that is \textit{regular} (uniformly
Lipschitz in its spatial variable).
We show that the regulariy of such decoupling field is closely related
to the bounded solution
to an associated \textit{characteristic BSDE}, a backward stochastic
Riccati-type equation with
superlinear growth in both components $Y$ and $Z$. We establish various
sufficient conditions
for the well-posedness of an ODE that dominates the characteristic
BSDE, which leads
to the existence of the desired regular decoupling random field, whence
the solvability of
the original FBSDE. A
synthetic analysis of the solvability is given, as a ``User's Guide,''
for a large class of FBSDEs
that are not covered by the existing methods. Some of them have
important implications in applications.
\end{abstract}

%
\end{frontmatter}

\section{Introduction}\label{sec1}

The theory of Backward Stochastic Differential Equations (BSDEs) and
Forward--Backward Stochastic Differential Equations (FBSDEs) have been
studied extensively for the past two decades, and its applications have
been found in many branches of
applied mathematics, especially the stochastic control theory and
mathematical finance. It has been noted, however, that while in many situations
the solvability of the original (applied) problems is essentially
equivalent to the
solvability of certain type of FBSDEs, these FBSDEs are often beyond
the scope
of any existing frameworks, especially when they are outside the
Markovian paradigm,
where the PDE tool becomes powerless. In fact,
the balance between the regularity of the coefficients and the time
duration, as well as
the nondegeneracy (of the forward diffusion), has been a longstanding
problem in the FBSDE literature,
especially in a general non-Markovian framework.
It has become increasingly clear that the theory now calls for new
insights and ideas that can lead to a better understanding of the problem
and hopefully to a unified solution scheme for the general FBSDEs.

A strongly coupled FBSDE takes the following form:
%
\begin{equation}\label{FBSDE}
\qquad\quad\cases{ \displaystyle X_t= x+\int_0^tb(s,X_s,
Y_s, Z_s)\,ds
\vspace*{3pt}\cr
\displaystyle\phantom{X_t=}{}
 +\int_0^t
\sigma(s,X_s,Y_s,Z_s)\,d B_s;
\vspace*{3pt}\cr
\displaystyle Y_t=g(X_T)+\int_t^Tf(s,X_s,
Y_s, Z_s) \,ds -\int_t^TZ_s\,dB_s,}\qquad
t\in[0,T],
\end{equation}
where $b$, $f$ and $\sigma $ are (progressively) measurable functions
defined on appropriate spaces, $B$ is a standard Brownian motion
and $g$ is a (possibly random) function that is defined on $\mathbb
{R}^n\times\Omega $
such that $g(x,\cdot )$ is $\cF_T$-measurable for each fixed $x$.

There have been three main methods to solve FBSDE (\ref{FBSDE}).
First, the
\textit{Method of Contraction Mapping}. This method, first used by
Antonelli \cite{fabio} and later detailed by Pardoux and Tang
\cite{PT}, works well when the duration $T$ is relatively small. Second,
the \textit{Four Step Scheme.} This was the first solution method that
removed restriction on the time duration for Markovian FBSDEs,
initiated by Ma, Protter and Yong \cite{mpy}. The trade-off is the
requirement on the regularity of the coefficients so that a
``decoupling'' quasi-linear PDE has a classical solution. Third, the
\textit{Method of Continuation.} This was a method that can treat
non-Markovian FBSDEs with arbitrary duration, initiated by Hu and Peng
\cite{hupeng} and Peng and Wu \cite{PW}, and later developed by Yong
\cite{yong1} and recently in \cite{yong4}. The main
assumption for this method is the so-called ``monotonicity conditions''
on the
coefficients, which is restrictive
in a different way. This method has been used widely in applications
(see, e.g.,
\cite{wu2,yu,wy}) because of its pure probabilistic nature.
We refer to the book of Ma and Yong \cite{mybk} for the detailed accounts
for all three methods. It is
worth noting that these three methods do not cover each other.

To make our motivation clearer, let us take a quick look at some main
difficulties in the FBSDE theory.
For example, consider the following simple FBSDE:
%
\begin{equation}
\label{eg1} X_t = x + \int_0^t
\sigma Z_s \,dW_s,\qquad Y_t = X_T
- \int_t^T Z_s \,dW_s.
\end{equation}
Clearly, the FBSDE has infinitely many solutions when $\sigma = 1$,
and is well-posed when $\sigma = 0$.
But more or less surprisingly, for $\sigma \neq0, 1$, none of the
three standard methods works for this FBSDE
when $T$ is arbitrarily large.
The FBSDE with such a feature has been encountered in many stochastic
control problems when diffusion
contains control, which is often the case in the optimal investment
problems in finance. Understanding its solvability is therefore
extremely desirable, especially when seeking the closed-loop optimal
control via Pontrygin's maximum
principle.
Another simple example, appeared in an earlier works of the fourth
author \cite{CZ} where the idea of \textit{method of optimal control}
(cf., e.g., \cite{mybk}) was adopted to study a Monte Carlo method for
FBSDEs, is of the following form:
%
\begin{eqnarray}
\label{eg2} X_t&=&x + \int_0^t
[a_sX_s + b_s Z_s ]\,ds + \int
_0^t \sigma _s \,dB_s;
\nonumber\\[-8pt]\\[-8pt]\nonumber
Y_t &=&h X_T + \int_t^T
[c_s X_s + d_s Z_s ]\,ds - \int
_t^T Z_s \,dB_s,
\end{eqnarray}
where $a$, $b$, $c$, $d$ and $\sigma $ are
stochastic processes, and $h$ is an $\cF_T$-random variable. Again,
this FBSDE is not covered by any existing method. However, as we
will see in Section~\ref{sec7} that the solvability of (\ref{eg1}) and (\ref{eg2}),
including an crucial estimate in \cite{CZ} regarding the solution to
(\ref{eg2}),
will be the easy consequences of our
general results. In fact, the work \cite{CZ} was the motivation for
\cite{Zhang}, which in turn motivated this paper.

The main goal of this paper is to
develop a strategy to construct a decoupling random field which will be
the key
to the solvability of general non-Markovian FBSDEs. Our starting point
is the work
of Delarue \cite{Delarue}, in a Markovian framework with $\sigma =
\sigma (t,x,y)$ being uniformly nondegenerate. In that case, an FBSDE
over arbitrary time duration
was solved under only Lipschitz conditions on the coefficients, by
combining nicely the Method of Contraction Mapping, the Four Step
Scheme, and some delicate PDE arguments. The idea was later extended by
Zhang \cite{Zhang} to the non-Markovian cases [again in the case
$\sigma = \sigma (t,x,y)$],
by using mainly probabilistic arguments, and with the help of some
compatibility conditions.
The main point is still, as in the Four Step Scheme, around finding a
function $u$ such that
%
\begin{equation}
\label{Yu} Y_t = u(t, X_t), \qquad t\in[0,T].
\end{equation}
Clearly, if the FBSDE (\ref{FBSDE}) is non-Markovian, then $u$ should
be a random field. The key issue here, as we shall argue, is the
existence of such a \textit{decoupling} random field that is uniformly
Lipschitz in its spatial variable. We will show that
the existence of such a random field is closely related to the
solvability of an associated BSDE (called the \textit{characteristic BSDE}
in this paper), and will
ultimately lead to the well-posedness of the original FBSDEs.
We shall provide a set of sufficient conditions for the existence of
such decoupling field, and show that most of the existing frameworks in
the literature could be
analyzed by using our criteria. Furthermore, we note that in the case
when the FBSDE is linear with
constant coefficients, some of our conditions are actually necessary.
In other words,
these conditions \textit{cannot} be improved.

A brief description of our plan is as follows.
Assume that the decoupling field $u$ exists and the FBSDE is
well-posed. Denote
$(X^x, Y^x, Z^x)$ to be the solution to FBSDE (\ref{FBSDE}) with
initial value $x$.
Then we argue that the derivative of $(X^x, Y^x, Z^x)$ with
respect to $x$, denoted by $(\nabla X, \nabla Y, \nabla Z)$, would
satisfy a linear ``variational FBSDE''
[see (\ref{variFBSDE}) below].
Since $Y^x_t = u(t, X^x_t)$ by (\ref{Yu}), we must have
$\nabla Y_t = u_x(t, X_t) \nabla X_t$, and thus $u_x(t, X_t) =
\nabla Y_t (\nabla X_t)^{-1}\stackrel{\triangle}{=}\hat Y_t $.
In other words, proving $u$ is uniformly Lipschitz continuous amounts
to finding
solutions to the linear FBSDE (\ref{variFBSDE}) such that $\hat Y$ is
uniformly bounded. Furthermore, one can check that $\hat Y$ actually satisfies
a BSDE [see (\ref{hatYBSDE}) below] which will be called the \textit{characteristic BSDE}
in this paper.
We note that this BSDE has superlinear growth in both components
of the solutions, thus it is itself a novel subject in BSDE theory, and
thus is interesting
in its own right.

Seeking the bounded solution to the characteristic BSDE over an arbitrary
time duration is by no means trivial, due to its superlinear growth behavior.
We shall accomplish this by studying two dominating ODEs [see (\ref
{ODE}) below], which bound $\hat Y$
from above and below, respectively. Although the ODEs also have the
combined complexity from its nonlinearity, superlinear growth, and the
singularity, it is much more tractable.
We shall give a set of
sufficient conditions to guarantee the existence of the solutions to
the ODEs, which in turn guarantees the
solvability of the original FBSDE (\ref{FBSDE}). Our results extend
those of
\cite{Zhang} in many ways, and we believe they are by far the most
general criteria for
the solvability of FBSDEs. As a byproduct, we also prove a comparison
theorem for the
decoupling random field over all time, thus confirming a common belief
(see, e.g.,
\cite{mybk,wu1,wx}).


There are several technical aspects in this paper that are worth
emphasizing. First,
unlike the linear FBSDEs studied in \cite{yong3}, where conditions
were made so that
the associated characteristic BSDE
is linear in $\hat Y$, or the so-called backward stochastic Riccati
equation, often seen in the linear-quadratic stochastic control
literature (see, e.g., \cite{KT} and \cite{Tang}) in which the growth
condition is quadratic in $\hat
Y$ but linear in $Z$, in the present case the generator has at least
quadratic growth on both
components.
To our best knowledge, such a case has not been investigated in the literature.
Second, our method requires the minimum assumptions on the
coefficients, and covers both Markovian and non-Markovian cases,
without having to go through the quasilinear PDEs and backward SPDEs
(see, e.g., \cite{Delarue,hmy,mpy,my,my1}). In an accompanying paper
\cite{MYZ}, however, we show that the
FBSDE has a uniformly Lipschitz continuous decoupling field (and thus
is well-posed)
if and only if the corresponding quasi-linear BSPDE has a uniformly
Lipschitz continuous Sobolev type weak solution. We hope that this
connection can enhance further
understanding on both FBSDEs and BSPDEs. Third, the method in this
paper is particularly effective for the cases where the forward
diffusion coefficient $\sigma $ depends on $Z$, which
has been avoided in many existing works, as it brings in some extra
complications for
the solvability analysis (see, e.g., \cite{Delarue,mybk}). Finally, in
this paper we
content ourselves for one-dimensional FBSDEs. In fact, the
characteristic BSDE becomes much more subtle in high-dimensional cases,
as it involves the combination of high-dimensional BSDEs with quadratic
growth (in $Z$) and high-dimensional backward stochastic Riccati
equations, each of which is very challenging. We hope to be able to
address this issue in our future publications.

The rest of the paper is organized as follows. In Section~\ref{sec2}, we
introduce the decoupling field and show how it leads to the
well-posedness of FBSDEs. In Section~\ref{sec3}, we heuristically discuss our
strategy for obtaining the uniformly Lipschitz continuity of the
decoupling field. In Section~\ref{sec4}, we study the relation between the
solvability of the linear variational FBSDE and its characteristic
BSDE, and in Section~\ref{sec5} we investigate the global solutions of the
dominant ODEs. In Section~\ref{sec6}, we investigate the well-posedness of
FBSDEs over small time duration, and in Section~\ref{sec7} we conclude our
well-posedness result for general FBSDEs over arbitrary time interval.
In Section~\ref{sec8}, we prove several further properties of FBSDEs. Finally,
in the \hyperref[app]{Appendix}, we complete some technical proofs.

\section{The decoupling field}\label{sec2}

Throughout this paper, we denote $(\Omega, \cF,\mathbb{P}; \mathbb
{F})$ to be a filtered
probability space on which is defined a Brownian motion $B=(B_t)_{t\ge
0}$. We assume
that $\mathbb{F}\stackrel{\triangle}{=}\mathbb{F}^B \stackrel
{\triangle}{=}\{\cF^B_t\}_{t \ge0}$, the natural filtration
generated by $B$,
augmented by the $\mathbb{P}$-null sets of $\cF$.
For any sub-$\sigma $-filed $\cG\subseteq\cF$, and $0\le p\le
\infty$, we denote
$L^p(\cG)$ to be the spaces of all $\cG$-measurable, $L^p$-integrable
random variables.
In what follows, we assume that
\textit{all processes involved are one-dimensional.}

Let $T>0$ be a fixed time horizon. We consider the general FBSDEs (\ref{FBSDE}), where the coefficients $b, \sigma, f, g$ are measurable
functions, and are
allowed to be random in general. For technical clarity, we shall make
use of the following \textit{standing assumptions} throughout the paper.

\begin{assum}
\label{assum-Lipschitz}
(i) The coefficients $b, \sigma, f\dvtx [0,T]\times\Omega \times\mathbb
{R}^3\mapsto\mathbb{R}$ are \mbox{$\mathbb{F}$-}progressively measurable,
for fixed $(x, y, z)\in\mathbb{R}^3$; and
the function $g\dvtx \mathbb{R}\times\Omega \mapsto\mathbb{R}$ is $\cF
_T$-measurable, for fixed $x\in\mathbb{R}$.
Moreover, the following integrability condition holds:
%
\begin{eqnarray}\label{I0}
I_0^2 &\stackrel{\triangle} {=}&\mathbb{E}
\biggl\{ \biggl(\int_0^T \bigl[\llvert b\rrvert
+ \llvert f\rrvert \bigr](t,0,0,0)\,dt \biggr)^2 +\int
_0^T \llvert \sigma \rrvert
^2(t,0,0,0)\,dt + \bigl\llvert g(0)\bigr\rrvert ^2 \biggr\}\hspace*{-30pt}
\nonumber\\[-8pt]\\[-8pt]\nonumber
&<& \infty.\hspace*{-30pt}
\end{eqnarray}

(ii) The coefficients $b,\sigma, f, g$ are
uniformly Lipschitz continuous in the spatial variable $(x, y, z)\in
\mathbb{R}^3$, uniformly in $\omega \in\Omega $, and with a common
Lipschitz constant $K_0>0$.
\end{assum}

To simplify notation, throughout the paper we denote
$\Theta \stackrel{\triangle}{=}(X, Y, Z)$.
Our purpose is to find $\mathbb{F}$-progressively measurable,
square-integrable processes $\Theta $, such that (\ref{FBSDE}) holds
for all $t\in[0,T]$, $\mathbb{P}$-a.s. However, to facilitate the
discussion, in what follows we often consider the FBSDE on a
subinterval $[t_1,t_2]$:
%
\begin{equation}
\label{FBSDEu} \cases{ \displaystyle X_t = \eta+ \int
_{t_1}^t b(s, \Theta _s) \,ds + \int
_{t_1}^t \sigma (s,\Theta _s)
\,dB_s;
\cr
\displaystyle Y_t = \varphi (X_{t_2})
+ \int_t^{t_2} f(s, \Theta _s) \,ds -
\int_t^{t_2} Z_s \,dB_s,}
\qquad t\in[t_1, t_2],
\end{equation}
where $\eta\in L^2(\cF_{t_1})$ and $\varphi (x,\cdot )\in L^2(\cF
_{t_2})$, for each fixed $x$.
We denote the solution to FBSDE (\ref{FBSDEu}), if exists, by $\Theta
^{t_1, t_2, \eta, \varphi }$.
In particular, we denote $\Theta ^{t,x}:= \Theta ^{t, T, x, g}$.

A well understood technique for solving an FBSDE, initiated in \cite
{mpy}, is to find
a ``decoupling function'' $u$ so that the solution $\Theta $ to the FBSDE
satisfies the relation~(\ref{eg1}). In Markovian cases, especially when
$\sigma = \sigma (t, x, y)$, it was shown that
$u$ is related to the solution to a quasilinear PDE, either in
classical sense
or in viscosity sense (cf., e.g., \cite{mpy,Delarue} or \cite{PT}). When the
coefficients are allowed to be random, special cases were also studied
and the function $u$ was found either as the solution to certain
backward stochastic PDEs (see, \cite{my,my1}), or as a random field
constructed by extending the localization technique of \cite{Delarue}
under certain compatibility conditions of the coefficients (see, \cite
{Zhang}). In the sequel, we call such random function $u$ the \textit{decoupling random field}
or simply \textit{decoupling field} of the FBSDE~(\ref{FBSDE}). More precisely, we have the following definition.

\begin{defn}
\label{defn-decoupling}
An $\mathbb{F}$-progressively measurable random field $u\dvtx  [0, T]\times
\mathbb{R}\times\Omega \mapsto\mathbb{R}$ with $u(T,x) = g(x)$ is
said to be a ``decoupling field'' of FBSDE (\ref{FBSDE}) if
%
%
%
there exists a constant $\delta >0$ such that, for any $0=t_1<t_2\le
T$ with $t_2-t_1\le\delta $ and any $\eta\in L^2(\cF_{t_1})$, the
FBSDE (\ref{FBSDEu}) with initial value $\eta$ and terminal condition
$u(t_2,\cdot )$ has a unique solution that satisfies (\ref{Yu}) for
$t\in[t_1, t_2]$, $\mathbb{P}$-a.s.

A decoupling field $u$ is called \textit{regular} if it is uniformly
Lipschitz continuous in~$x$.
\end{defn}

By a slight abuse of notation, we shall denote the solution in
Definition \ref{defn-decoupling} by
$\Theta ^{t_1, t_2, \eta, u}$. One should note that the existence of
the (regular) decoupling field implies the well-posedness of the FBSDE
over a small interval, which is usually guaranteed by the Method of
Contraction Mapping given the Assumption \ref{assum-Lipschitz}.
The following result shows the significance of the existence of the
decoupling field for the well-posedness for FBSDEs over an arbitrary duration.

%
\begin{teo}
\label{teo-decoupling}
Assume that Assumption \ref{assum-Lipschitz} holds, and that there
exists a decoupling field $u$ for FBSDE (\ref{FBSDE}). Then FBSDE
(\ref{FBSDE}) has a unique solution $\Theta $ and (\ref{Yu}) holds
over an arbitrary duration $[0, T]$.
\end{teo}

\begin{pf}
Let $T>0$ be given. Consider a partition: $0=t_0<\cdots<t_n=T$ of $[0,T]$ such that $t_{i+1}-t_{i}\le\delta $, $i=0,\ldots, n-1$, where $\delta $ is the constant in Definition~\ref{defn-decoupling}.

Define $X_{t_0}\stackrel{\triangle}{=}x$, $Y_{t_0}\stackrel
{\triangle}{=}u(0,x)$, and for $i=0,\ldots, n-1$, define recursively
\begin{eqnarray*}
\Theta _t\stackrel{\triangle} {=}\Theta ^{t_{i}, t_{i+1}, X_{t_{i}},
u}_t,
\qquad t\in(t_{i}, t_{i+1}].
\end{eqnarray*}
Then $\Theta $ would solve FBSDE (\ref{FBSDE}) if they could be
``patched'' together. But note that
\begin{eqnarray*}
X_{t_{i}+} &=& X^{t_{i}, t_{i+1}, X_{t_{i}}, u}_{t_{i}+} = X^{t_{i},
t_{i+1}, X_{t_{i}}, u}_{t_{i}}
= X_{t_i};
\\
Y_{t_{i}+} &=& Y^{t_{i}, t_{i+1}, X_{t_{i}}, u}_{t_{i}+} = Y^{t_{i},
t_{i+1}, X_{t_{i}}, u}_{t_{i}}
= u(t_i, X_{t_i}) = Y^{t_{i-1}, t_{i},
X_{t_{i-1}}, u}_{t_{i}} =
Y_{t_i}.
\end{eqnarray*}
That is, $(X, Y)$ is continuous on $[0, T]$. Moreover, $u(T, x) =
g(x)$, then by (\ref{FBSDEu}) one can check straightforwardly that
$\Theta $ satisfies FBSDE (\ref{FBSDE}) on $[0, T]$, proving the
existence. Furthermore, from our construction it is clear that (\ref
{Yu}) holds.

We now prove the uniqueness. Let $\tilde\Theta $ be an arbitrary
solution to FBSDE (\ref{FBSDE}). Note that $\tilde\Theta $ satisfies
FBSDE (\ref{FBSDEu}) on $[t_{n-1}, t_n]$ with initial condition
$\tilde X_{t_{n-1}}$. Then by the definition of the decoupling field,
we have $\tilde Y_{t_{n-1}} = u(t_{n-1}, \tilde X_{t_{n-1}})$. This
implies that $\tilde\Theta $ satisfies FBSDE (\ref{FBSDEu}) on
$[t_{n-2}, t_{n-1}]$ with initial condition $\tilde X_{t_{n-2}}$. Then
we have $\tilde Y_{t_{n-2}} = u(t_{n-2}, \tilde X_{t_{n-2}})$.
Repeating the arguments backwardly in time, we obtain that $\tilde
Y_{t_i} = u(t_i,\tilde X_{t_i})$, $i=n,\ldots, 0$. Now consider
FBSDE (\ref{FBSDEu}) on $[t_0, t_1]$. Since $\tilde X_{t_0} = x =
X_{t_0}$, by the uniqueness of solutions we know that $\tilde\Theta =
\Theta $ on $[t_0, t_1]$. In particular, $\tilde X_{t_1} = X_{t_1}$,
and thus the corresponding FBSDEs (\ref{FBSDEu}) on $[t_1, t_2]$ have
the same initial condition. Repeating the arguments, this time
forwardly for $i=1,\ldots, n$, we see that $\tilde\Theta = \Theta
$ on $[0, T]$, and thus the solution is unique.
\end{pf}

We conclude this section by making the following observations.

\begin{rem}\label{rem-decoupling}
(i) Definition \ref{defn-decoupling} and Theorem \ref{teo-decoupling} can be extended to higher-dimensional cases (with the
constant $\delta $ possibly depending on the dimensions as well), and
the proof stays exactly the same.

(ii) By the uniqueness in Theorem~\ref{teo-decoupling}, it is obvious
that the decoupling field, if exists, is also unique. In fact, it is
clear that $u(t,x) = Y^{t, x}_t$.
\end{rem}

\begin{rem}
\label{rem-Lipschitz}
A typical condition for well-posedness of FBSDEs over small time
interval is the uniform Lipschitz continuity of the terminal condition.
Therefore, the main goal of this paper is to provide sufficient
conditions which guarantee the existence of the regular decoupling
field $u$.
Such a feature was also observed from a different angle in \cite{MYZ},
in which we characterize the regular decoupling field $u$ as a Sobolev
type weak solution to certain \textit{backward stochastic PDE} that is
Lipschitz in $x$. We note that the idea of ``decoupling device'' was
also used for
linear FBSDEs in \cite{yong3}. But in that work the
uniform Lipschitz continuity was not studied.
\end{rem}

\section{Some Heuristic analysis}\label{sec3}\label{sect-heuristic}

From Theorem~\ref{teo-decoupling} and Remark \ref{rem-Lipschitz}, it
is easy to see that the issue of the well-posedness of FBSDE (\ref
{FBSDE}) can be decomposed into two parts. First, the well-posedness on
small time interval, and second, finding a decoupling field $u$ that is
uniformly Lipschitz continuous in its spatial variable. The
first issue was more or less ``classical'' (see, e.g., \cite{fabio}),
but we will fine-tune it in Section~\ref{sect-smallT} to suit our
purpose. The second issue, however, is much more subtle, and is the
main focus of this paper. In this section, we first give a heuristic
analysis, from which several fundamental problems will be formulated,
and their proofs will be carried out in Sections~\ref{sect-linear} and
\ref{sect-ODE} below. A synthetic analysis will then be given in
Section~\ref{sect-summary}.

We first introduce some notation: for $\theta _j:= (x_j, y_j, z_j)$,
$j=1,2$, and for $\varphi = b, \sigma, f$, denote
%
\begin{eqnarray}\label{tildephi}
\tilde h(x_1, x_2) &\stackrel{\triangle}{=}& \bigl[g(x_1)-g(x_2)\bigr]\slash[x_1-x_2];\nonumber
\\
\tilde\varphi _1(t, \theta _1, \theta _2)&\stackrel {\triangle} {=}& \bigl[\varphi (t,x_1, y_1,
z_1) - \varphi (t, x_2, y_1,z_1)\bigr]\slash[x_1-x_2];
\nonumber\\[-8pt]\\[-8pt]\nonumber
\tilde\varphi _2(t, \theta _1, \theta _2)&\stackrel {\triangle} {=}& \bigl[\varphi (t,x_2, y_1,
z_1) - \varphi (t, x_2, y_2,z_1)\bigr]\slash[y_1-y_2];
\\
\tilde\varphi _3(t, \theta _1, \theta _2)&\stackrel {\triangle} {=}& \bigl[\varphi (t,x_2, y_2,
z_1) - \varphi (t, x_2, y_2,
z_2)\bigr]\slash[z_1-z_2].\nonumber
\end{eqnarray}
Here and in the sequel, for any Lipschitz continuous function $\varphi
(x)$, when $x_1=x_2$ we will always take the convention that
%
\begin{equation}
\label{00} { \varphi (x) - \varphi (x) \over x-x}:= \liminf_{\tilde x\to x}
{\varphi (\tilde x) - \varphi (x)\over\tilde x-x}.
\end{equation}

Our main idea to decouple the FBSDE (\ref{FBSDE}) is as follows.
Assume that there exists a decoupling field $u=u(t,x)$ that is
uniformly Lipschitz continuous in $x$ and (\ref{Yu}) holds. Assume
also that (\ref{FBSDE}) is well-posed on $[0,T]$, with $X_0=x$ for any
$x$. Given $x_i$, $i=1, 2$, let $\Theta ^i$ denote the unique solution
to (\ref{FBSDE}) with initial condition $x_i$, and. By slightly
%
\begin{equation}
\label{tdX} \nabla \Theta \stackrel{\triangle} {=}{\Theta ^1-\Theta ^2 \over
x_1 - x_2},\qquad
\nabla u(t) \stackrel{\triangle} {=}{u(t, X^1_t)- u(t,
X^2_t)\over X^1_t- X^2_t}.
\end{equation}
Since $Y^i_t = u(t, X^i_t)$, $i=1, 2$, one must have
%
\begin{equation}
\label{uxXY} \nabla Y_t= \nabla u(t) \nabla X_t,
\end{equation}
and one can check immediately that $\nabla \Theta $ satisfies the
following ``variational \mbox{FBSDE}:''
%
\begin{equation}
\label{variFBSDE}
\qquad\cases{ \displaystyle \nabla X_t= 1+\int
_0^t (b_1\nabla
X_s+b_2\nabla Y_s+b_3\nabla
Z_s)\,ds
\vspace*{5pt}\cr
\displaystyle\phantom{\nabla X_t=} {} +\int_0^t (\sigma
_1\nabla X_s+\sigma _2\nabla
Y_s+\sigma _3\nabla Z_s)\,dB_s;
\vspace*{5pt}\cr
\displaystyle \nabla Y_t= h\nabla X_T +\int
_t^T (f_1\nabla
X_s+f_2\nabla Y_s+f_3\nabla
Z_s)\,ds
\vspace*{5pt}\cr
\displaystyle\phantom{\nabla Y_t=}
{}  -\int_t^T \nabla
Z_s\,dB_s,}\qquad
t\in[0,T],
\end{equation}
where
$h \stackrel{\triangle}{=}\tilde h(X^1_T, X^2_T)$ and $\varphi _i(t)
\stackrel{\triangle}{=}\tilde\varphi _i(t, \Theta ^1_t, \Theta
^2_t)$, $i=1,2,3$, $ \varphi = b$, $\sigma $, $f$, respectively.
We note here that $b_i, \sigma _i, f_i$, $i=1,2,3$, are $\mathbb
{F}$-adapted processes and $h$ is a $\cF_T$-measurable random
variable, and they are all bounded, thanks to Assumption~\ref{assum-Lipschitz}.

Furthermore, in light of (\ref{uxXY}) we see that a decoupling field
$u$ being regular (i.e., uniformly Lipschitz continuous in $x$) is
essentially equivalent to $\hat Y_t \stackrel{\triangle}{=}\nabla
Y(\nabla X)^{-1}$ being uniformly bounded. Thus, let us assume $\nabla
X\neq0$ and denote
%
\begin{eqnarray}\label{hatYZ}
\hat Y_t &\stackrel{\triangle} {=}& \nabla Y_t/\nabla X_t\quad\mbox{and}
\nonumber\\[-8pt]\\[-8pt]\nonumber
\hat Z_t &\stackrel{\triangle} {=}& \bigl[\nabla Z_t - \hat Y_t(\sigma
_1\nabla X_t + \sigma _2 \nabla Y_t + \sigma _3\nabla Z_t)\bigr]\slash \nabla
X_t,
\end{eqnarray}
or equivalently,
%
\begin{equation}
\label{YZ} \nabla Y_t = \hat Y_t\nabla
X_t,\qquad \nabla Z_t = {\hat Z_t + \hat
Y_t(\sigma _1 + \sigma _2 \hat Y_t)\over1-\sigma _3\hat Y_t}\nabla
X_t.
\end{equation}
%
A simple application of It\^o's formula to $\hat Y_t$, assuming
$\sigma _3 \hat Y \neq1$, yields that
%
\begin{equation}
\label{hatYBSDE} \quad\hat Y_t = h + \int_t^T
\bigl[F_s(\hat Y_s) + G_s(\hat
Y_s) \hat Z_s + \Lambda _s(\hat
Y_s)\llvert \hat Z_s\rrvert ^2 \bigr]\,ds -
\int_t^T \hat Z_s
\,dB_s,
\end{equation}
where
%
\begin{eqnarray}
\label{F} F_s(y) &\stackrel{\triangle} {=}&
f_1+f_2y + y(b_1+b_2y)+
{(f_3 + b_3y)
y(\sigma _1 + \sigma _2 y)\over1-\sigma _3y};
\nonumber
\\
G_s(y) &\stackrel{\triangle} {=}& \sigma _1+\sigma
_2y + {f_3 + b_3y
+ \sigma _3 y(\sigma _1 + \sigma _2 y)
\over 1-\sigma _3y}
\nonumber\\[-8pt]\\[-8pt]\nonumber
&=& {(\sigma _1+f_3) +(\sigma _2+ b_3)y \over
1-\sigma _3y};
\\
\Lambda _s(y) &\stackrel{\triangle} {=}& {\sigma _3 \over 1-\sigma
_3y}.
\nonumber
\end{eqnarray}
Equation (\ref{hatYBSDE}) is clearly a legitimate BSDE, even
without assuming $\nabla X\neq0$. We shall call this BSDE the ``\textit{Characteristic BSDE}'' of the linear variational FBSDE~(\ref{variFBSDE})
[or of the original FBSDE (\ref{FBSDE})], and their connection
will be studied rigorously in the next section.
We note that the identities in (\ref{YZ}) and the desired Lipschitz
property of the decoupling field $u$ tell us that we should look for
conditions under which the BSDE (\ref{hatYBSDE}) has a solution $(\hat
Y, \hat Z)$ such that
%
\begin{eqnarray}
\label{hatYproperty} \mbox{both $\hat Y$ and $(1-\sigma _3 \hat
Y)^{-1}$ are bounded.}
\end{eqnarray}

\begin{rem}
It is worth noting that the BSDE (\ref{hatYBSDE}) is nonstandard in
several aspects. Most notable is that its generator has at least
quadratic growth in both $Y$ and $Z$, thus it can be thought of as a
Backward Stochastic Riccati Equations (BSRE) with quadratic growth in
$Z$, which to our best knowledge, has not been studied in literature.

Besides the commonly cited reference of BSDEs with quadratic growth in
$Z$ (e.g., \cite{Kobylanski,BH}), the following special cases
of (\ref{hatYBSDE}) are worth mentioning. In \cite{Tang}, the BSRE
with linear growth in $Z$ was studied in the context of stochastic LQ
(linear-quadratic) problem, in which the FBSDE is a natural consequence
of the stochastic maximum principle. The characteristic BSDE (\ref
{hatYBSDE}) was also observed in \cite{yong3}, where the linear FBSDEs
were considered. But some special assumptions were made so that the
BSDE has linear growth in $Y$. Finally, in \cite{Zhang} certain
compatibility conditions were also added so that (\ref{hatYBSDE})
becomes a standard BSDE, and thus its well-posedness was not an issue.
Our results will
contain those of \cite{Tang,yong3} and \cite{Zhang} as
special cases.
\end{rem}

We conclude this section by outlining the strategy for obtaining the
{a priori} uniform estimate of $\hat Y$, which is crucial for
finding the solution of (\ref{hatYBSDE}) satisfying~(\ref
{hatYproperty}). To begin with, for any bounded random variable $\xi$,
define its deterministic
upper and lower bounds by
%
\begin{eqnarray}
\label{xibound}
\overline\xi &\stackrel{\triangle} {=}&\esssup\xi\stackrel{\triangle
} {=}\inf\{a\in\mathbb{R}\dvtx  \xi\le a, \mbox{ a.s.}\},
\nonumber\\[-8pt]\\[-8pt]\nonumber
\underline\xi &\stackrel{\triangle} {=}&\essinf\xi\stackrel {\triangle} {=}\sup\{a\in\mathbb{R}\dvtx
\xi\ge a,\mbox{ a.s.}\}.
\end{eqnarray}
For any $\theta _j\stackrel{\triangle}{=}(x_j, y_j, z_j)$, $j=1,2$,
we define $F(\theta _1,\theta _2; t,y)$ by replacing the coefficients
$\varphi _i$ in (\ref{F}) with $\tilde\varphi _i(t, \theta
_1,\theta _2)$ defined in (\ref{tildephi}), $i=1,2,3$, $\varphi
=b,\sigma, f$. We then define
%
\begin{eqnarray}
\label{Fbar} \overline h &\stackrel{\triangle} {=}&\esssup \Bigl(\sup
_{x_1 \neq
x_2} \tilde h(x_1,x_2) \Bigr),\nonumber
\\
\underline h &\stackrel{\triangle} {=}& \essinf \Bigl(\inf_{x_1 \neq
x_2} \tilde h(x_1, x_2) \Bigr),
\nonumber\\[-8pt]\\[-8pt]\nonumber
\overline F(t,y) &\stackrel{\triangle} {=}&\esssup \Bigl(\sup
_{x_1
\neq x_2, y_1\neq y_2, z_1\neq z_2} F(\theta _1,\theta _2; t,y)
\Bigr),
\\
\underline F(t,y) &\stackrel{\triangle} {=}&\essinf \Bigl(\inf
_{x_1
\neq x_2, y_1\neq y_2, z_1\neq z_2} F(\theta _1, \theta _2; t,y)
\Bigr).
\nonumber
\end{eqnarray}
Here, we should remark that
$\overline F(t,y)$ is a \textit{deterministic} function, and we should
note its notational
difference from the possibly random processes, for example, $F_t(y)$,
$G_t(y)$, etc., appeared previously.
We\vspace*{2pt} have the following a priori estimate of~$\hat Y$.

\begin{lem}
\label{lem-hatYbound} Let Assumption \ref{assum-Lipschitz} hold.
Assume that the BSDE (\ref{hatYBSDE}) has a solution $(\hat Y, \hat
Z)$, and the following ordinary differential equations (ODEs) admit
solutions $\overline\by, \underline\by$:
%
\begin{eqnarray}
\label{ODE} \overline\by_t = \overline h + \int
_t^T \overline F(s, \overline\by
_s) \,ds,\qquad \underline\by_t = \underline h + \int
_t^T \underline F(s, \underline
\by_s) \,ds.
\end{eqnarray}
Assume further that $\hat Y$, $\overline\by$ and $\underline\by$
all satisfy (\ref{hatYproperty}).
Then $\underline\by_t \le\hat Y_t \le\overline\by_t$, for all
$t\in[0,T]$, $\mathbb{P}$-a.s.
\end{lem}

\begin{pf}
Denote
$
\tilde G_t(z) \stackrel{\triangle}{=}G_t(\hat Y_t) z + \Lambda
_t(\hat Y_t) z^2$. Note that $(\hat Y, \hat Z)$ satisfies
the following BSDE:
\begin{eqnarray*}
Y_t &=& h + \int_t^T
\bigl[F_s(Y_s) + \tilde G_s(Z_s)
\bigr]\,ds - \int_t^T Z_s
\,dB_s
\end{eqnarray*}
and $(\overline\by, 0)$ satisfy the following BSDE:
\begin{eqnarray*}
Y_t &=& \overline h + \int_t^T
\bigl[\overline F(s, Y_s) + \tilde G_s(Z_s)
\bigr]\,ds - \int_t^T Z_s
\,dB_s.
\end{eqnarray*}
Let $C>0$ be the common upbound of $\llvert  \hat Y\rrvert  $, $\llvert  1-\sigma _3\hat
Y\rrvert  ^{-1}$, $\llvert  \overline\by\rrvert  $,
$\llvert  1-\sigma _3\overline\by\rrvert  ^{-1}$, $\llvert  \underline\by\rrvert  $, and
$\llvert  1-\sigma _3\underline\by\rrvert  ^{-1}$.
Note that $F$ is uniformly Lipschitz continuous in $y$ in the set $\{y\dvtx
\llvert  y\rrvert  \le C, \llvert  1-\sigma _3 y\rrvert  ^{-1}\le C\}$. It then follows from the
comparison theorem for quadratic BSDEs (see, e.g., \cite{Kobylanski})
that $\hat Y \le
\overline\by$. Similarly, we have $\hat Y \ge\underline\by$.
\end{pf}

Combining the discussions in Sections~\ref{sec2} and \ref{sec3}, especially Lemma~\ref{lem-hatYbound}, it is now clear
that finding the uniform Lipschitz decoupling random field $u$ will
eventually come down to finding conditions so that the ODEs in (\ref
{ODE}) admit nonexplosive solutions over the arbitrarily
prescribed duration $[0,T]$. In the rest of the paper,
we shall call the ODEs in (\ref{ODE}) the ``\textit{dominating ODEs}'' of
BSDE (\ref{hatYBSDE}), whose well-posedness will be the main subject
of Section~\ref{sect-ODE}.

\section{The characteristic BSDE}\label{sec4}\label{sect-linear}

In this section, we study the connection between well-posedness of
the linear variational FBSDE (\ref{variFBSDE}) and the corresponding
characteristic BSDE (\ref{hatYBSDE}). Such a relation is not only
interesting in its own right, but also important for us to
construct the desired regular decoupling field in Sections~\ref{sec6}~and~\ref{sec7} below. We should note that the variational
FBSDE (\ref{variFBSDE}) coincides with the original FBSDE if (\ref
{FBSDE}) is actually linear.

For notational simplicity, we denote $(\cX, \cY, \cZ):= (\nabla X,
\nabla Y, \nabla Z)$ and then the variational FBSDE (\ref
{variFBSDE}) becomes the following linear FBSDE
with random coefficients:
%
\begin{equation}
\label{linearFBSDE} \cases{ \displaystyle \cX_t= 1+\int
_0^t (b_1\cX_s+b_2
\cY_s+b_3\cZ _s)\,ds
\vspace*{3pt}\cr
\displaystyle \phantom{\cX_t=}{} +\int_0^t
(\sigma _1\cX_s+\sigma _2
\cY_s+\sigma _3\cZ_s)\,dB_s;
\vspace*{3pt}\cr
\displaystyle \cY_t= h\cX_T +\int_t^T
(f_1\cX_s+f_2\cY_s+f_3
\cZ_s)\,ds -\int_t^T
\cZ_s\,dB_s.}
\end{equation}
In this case, (\ref{hatYZ}) and (\ref{YZ}) become
%
\begin{eqnarray}
\label{hatYZ-linear} \hat Y &\stackrel{\triangle} {=}&\cY\slash\cX,\qquad \hat Z
\stackrel{\triangle} {=}\bigl[\cZ- \hat Y(\sigma _1\cX+ \sigma
_2 \cY + \sigma _3 \cZ)\bigr]\slash\cX,
\\
\label{YZ-linear} \cY&=& \hat Y \cX,\qquad \cZ= \bigl[\hat Z + \hat Y(\sigma
_1 + \sigma _2 \hat Y) \bigr] \cX\slash[1-\sigma
_3\hat Y].
\end{eqnarray}

The original Assumption \ref{assum-Lipschitz} can be translated into
the following assumption.

\begin{assum}
\label{assum-linear}
Assume $b_i, \sigma _i, f_i$, $i=1,2,3$, are $\mathbb{F}$-adapted
processes, $h$~is a $\cF_T$-measurable random variable, and they are
all bounded.
\end{assum}

The following spaces are important in our discussion. For $p\ge1$, denote
%
\begin{eqnarray}\label{space}
\mathbb{L}^p &\stackrel{\triangle} {=}& \biggl\{\Theta\dvtx  \llVert \Theta \rrVert _{\mathbb{L}^p}^p\stackrel{\triangle} {=}
\mathbb{E} \biggl\{\sup_{0\le t\le T}\bigl[\llvert X_t
\rrvert ^p + \llvert Y_t\rrvert ^p\bigr] +
\biggl(\int_0^T\llvert Z_t\rrvert
^2 \,dt \biggr)^{p/2} \biggr\} <\infty \biggr\};\hspace*{-30pt}
\nonumber\\[-8pt]\\[-8pt]\nonumber
\widehat{\mathbb{L}}_p &\stackrel {\triangle} {=}&\bigcup _{q>p}\mathbb{L}^q.
\end{eqnarray}

We begin our discussion with the following observation. For any
$\mathbb{F}$-adapted process $u$
such that $\int_0^T \llvert  u_t\rrvert  ^2\,dt<\infty$, $\mathbb{P}$-a.s., we define
%
\begin{eqnarray}
\label{Mu} M^u_t \stackrel{\triangle} {=}\exp
\biggl\{\int_0^t u_s
\,dB_s - {1\over
2}\int_0^t
\llvert u_s\rrvert ^2 \,ds \biggr\}.
\end{eqnarray}
Consider the following simplified form of (\ref{hatYBSDE}):
%
\begin{eqnarray}
\label{BMO-BSDE} \hat Y_t = h + \int_t^T
\bigl[\alpha _s + \beta_s \hat Y_s +
\gamma _s \hat Z_s + \lambda _s \llvert
\hat Z_s\rrvert ^2\bigr]\,ds - \int_t^T
\hat Z_s \,dB_s,
\nonumber\\[-8pt]\\[-8pt]
\eqntext{t\in[0,T],}
\end{eqnarray}
%
where $\alpha, \beta, \gamma, \lambda $ are $\mathbb{F}$-adapted
processes and $h$ is an $\cF_T$-measurable random variable, all
bounded. Then it is well known (see, e.g., \cite{BH}) that the BSDE
(\ref{BMO-BSDE}) admits a unique solution $(\hat Y, \hat Z)$ such
that, for some constant $C$ depending on the bounds of $\alpha, \beta, \gamma, \lambda, h$ and $T$, 
%
\begin{eqnarray*}
\llvert \hat Y_t\rrvert \le C 
\quad\mbox{and}\quad
\mathbb{E}_t \biggl\{\int_t^T
\llvert \hat Z_s\rrvert ^2 \,ds \biggr\} \le C.
\end{eqnarray*}
%
Furthermore, applying some BMO analysis (cf. \cite{KS}, Lemma~4 and Theorem~1), one shows that
there exists a constant $\varepsilon >0$, depending also on the bounds
of the coefficients and $T$, such that
%
\begin{eqnarray}
\label{BMO} \mathbb{E} \biggl\{\exp \biggl(\varepsilon \int
_0^T \llvert \hat Z_t\rrvert
^2\,dt \biggr) +\bigl\llvert M^{\lambda \hat Z}_T\bigr\rrvert
^{1+\varepsilon } \biggr\}<\infty.
\end{eqnarray}
Consequently, $M^{\lambda \hat Z}$ is a true martingale.


Bearing this observation in mind, we now give the main result of this section.

\begin{teo}
\label{teo-linear-wellposed} Assume Assumption \ref{assum-linear}
holds.

\textup{(i)} If\vspace*{2pt} the BSDE (\ref{hatYBSDE}) has a solution $(\hat Y, \hat Z)$
such that (\ref{hatYproperty}) holds,
then the FBSDE (\ref{linearFBSDE}) has a solution $(\cX, \cY, \cZ
)\in\widehat{\mathbb{L}}_1$ such that $\cX\neq0$ and (\ref{YZ-linear})
holds.\vspace*{2pt}

\textup{(ii)} Conversely, if the FBSDE (\ref{linearFBSDE}) has a solution $(\cX, \cY, \cZ)\in\widehat{\mathbb{L}}_1$ such that
%
\begin{eqnarray}
\label{XYproperty} \llvert \cY_t\rrvert \le C\llvert \cX_t
\rrvert, \qquad \llvert \cX_t\rrvert \le C\llvert \cX_t -
\sigma _3 \cY_t\rrvert,
\end{eqnarray}
then $\cX\neq0$, and the processes $(\hat Y, \hat Z)$ defined by
(\ref{hatYZ-linear}) satisfies BSDE (\ref{hatYBSDE}) and~(\ref{hatYproperty}).
\end{teo}

\begin{pf}
(i) In light
of (\ref{YZ-linear}), we consider the following SDE:
%
\begin{eqnarray}
\label{dX} d \cX_t &=& \cX_t \biggl[b_1
+ b_2 \hat Y_t + b_3 {\hat Z_t + \hat
Y_t(\sigma _1 + \sigma _2 \hat Y_t)\over1-\sigma _3\hat Y_t}
\biggr]\,dt \nonumber
\\
&&{} + \cX_t \biggl[\sigma _1 + \sigma
_2 \hat Y_t + \sigma _3 {\hat Z_t +
\hat Y_t(\sigma _1 + \sigma _2 \hat Y_t)\over1-\sigma _3\hat
Y_t}
\biggr]\,dB_t
\\
&=&\cX_t \bigl\{H_t(\hat Y_t,\hat
Z_t)\,dt + \bigl[I_t(\hat Y_t) + \Lambda
_t(\hat Y_t) \hat Z_t\bigr]\,dB_t
\bigr\},\nonumber
\end{eqnarray}
where
\begin{eqnarray*}
H_t(y,z) \stackrel{\triangle} {=} \biggl[b_1 +
b_2 y + b_3 {z +
y(\sigma _1 +
\sigma _2 y)\over1-\sigma _3y} \biggr]; \qquad
I_t(y) \stackrel{\triangle} {=}{\sigma _1 +
\sigma _2y\over1-\sigma _3y}.
\end{eqnarray*}
It is then easy to check that
%
\begin{eqnarray}
\label{X} \cX_t &=& \exp \biggl\{\int_0^t
\bigl[I_s(\hat Y_s) + \Lambda _s(\hat
Y_s) \hat Z_s\bigr]\,dB_s\nonumber
\\
&&\hspace*{20pt} {}
+ \int
_0^t \biggl[H_s(\hat
Y_s,\hat Z_s) - {1\over2}
\bigl[I_s(\hat Y_s)+ \Lambda _s(\hat
Y_s) \hat Z_s \bigr]^2 \biggr]\,ds \biggr\}
\\
&=& M^{\Lambda (\hat Y) \hat Z}_t M_t^{I(\hat Y)}\exp \biggl\{
\int_0^t \bigl[H_s(\hat
Y_s, \hat Z_s)-I_s(\hat Y_s)
\Lambda _s(\hat Y_s) \hat Z_s\bigr]\,ds
\biggr\}.\nonumber
\end{eqnarray}
Clearly, $\cX>0$. Furthermore, since (\ref{hatYproperty}) implies
that in (\ref{dX})
$\Lambda (\hat Y)$, $I(\hat Y)$ are bounded and $H(\hat Y, \hat Z)$
has a linear growth in
$\hat Z$, and
(\ref{BMO}) implies
%
\begin{eqnarray}
\label{hatZp} \mathbb{E} \biggl\{\sup_{0\le t\le T} \bigl\llvert
M_t^{I(\hat Y)}\bigr\rrvert ^p + \exp \biggl(p\int
_0^T \bigl[1 + \llvert \hat Z_t
\rrvert \bigr]\,dt \biggr) \biggr\}<\infty
\nonumber\\[-8pt]\\[-8pt]
\eqntext{\mbox{for any }p>1,}
\end{eqnarray}
we deduce from (\ref{X}) that, for $\varepsilon $ in (\ref{BMO})
[noting that $({2(1+\varepsilon )\over2+\varepsilon },
{2(1+\varepsilon )\over\varepsilon })$ are conjugates],
%
\begin{eqnarray}\label{X1eps}
&& \mathbb{E} \Bigl\{\sup_{0\le t\le T}\llvert
\cX_t\rrvert ^{1+{\varepsilon/2}} \Bigr\}\nonumber
\\
&&\qquad \le \Bigl(E \Bigl\{\sup
_{0\le t\le T}\bigl\llvert M_t^{\Lambda (\hat
Y) \hat Z}\bigr\rrvert^{1+\varepsilon } \Bigr\} \Bigr)^{(2+\varepsilon)/(2(1+\varepsilon ))}
\nonumber\\[-8pt]\\[-8pt]\nonumber
&&\quad\qquad{}  \times
\Bigl(E \Bigl\{\sup_{0\le t\le T}\bigl\llvert M_t^{I(\hat Y)}
\bigr\rrvert ^{(2+\varepsilon )(1+\varepsilon )/\varepsilon } e^{{C(2+\varepsilon )(1+\varepsilon )/\varepsilon }\int_0^T [1
+ \llvert  \hat Z_t\rrvert  ]\,dt} \Bigr\} \Bigr)^{\varepsilon /(2(1+\varepsilon))}\hspace*{-15pt}
\\
&&\qquad < \infty.\nonumber
\end{eqnarray}
Now if we define $(\cY, \cZ)$ by (\ref{YZ-linear}), then $(\cX, \cY, \cZ)$ satisfy
(\ref{linearFBSDE}) and, by (\ref{BMO}) again,
%
\begin{equation}
\label{YZ1eps} \mathbb{E} \biggl\{\sup_{0\le t\le T}\llvert
\cY_t\rrvert ^{1+{\varepsilon/2}} + \biggl(\int_0^T
\llvert \cZ_t\rrvert ^2\,dt \biggr)^{1+{\varepsilon /4}} \biggr\}
< \infty.
\end{equation}
That is, $(\cX, \cY, \cZ)\in\mathbb{L}^{1+{\varepsilon /4}}
\subset\widehat{\mathbb{L}}_1$, proving (i).

(ii) We now assume that FBSDE (\ref{linearFBSDE}) has a solution $(\cX, \cY, \cZ)\in\widehat{\mathbb{L}}_1$ such that~(\ref{XYproperty}) holds. Denote $\tau _n \stackrel{\triangle
}{=}\inf\{t\dvtx  \cX_t = {1\over n}\}\wedge T$, $\tau \stackrel
{\triangle}{=}
\inf\{t\dvtx  \cX_t = 0\}\wedge T$, and define $\hat Y, \hat Z$ by (\ref
{hatYZ-linear}). Clearly, the assumption (\ref{XYproperty}) implies
that $\hat Y$ satisfies (\ref{hatYproperty}) in~$[0,\tau )$, and
applying It\^{o}'s
formula we see that $(\hat Y, \hat Z)$ satisfies
\begin{eqnarray*}
d\hat Y_t = - \bigl[F_t(\hat Y_t) +
G_t(\hat Y_t) \hat Z_t + \Lambda
_t(\hat Y_t) \llvert \hat Z_t\rrvert
^2 \bigr] \,dt + \hat Z_t \,dB_t,\qquad t\in [0,\tau
).
\end{eqnarray*}
Note that the boundedness of $\hat Y$ implies that the above SDE is
actually of the form of (\ref{BMO-BSDE}), and at least on $[0,\tau
_n)$ the stochastic integral $\int_0^\cdot \hat Z_s \,dB_s$ is a
true martingale. Thus, we can apply the same argument there to obtain
the bound (\ref{BMO}) on
$[0,\tau _n)$:
%
\begin{eqnarray*}
\mathbb{E} \biggl\{\exp \biggl(\varepsilon \int_0^{\tau _n}
\llvert \hat Z_t\rrvert ^2\,dt \biggr) + \bigl\llvert
M^{\Lambda (\hat Y)\hat Z}_{\tau
_n}\bigr\rrvert ^{1+\varepsilon } \biggr\} \le C<
\infty.
\end{eqnarray*}
Note that the constants $\varepsilon $ and $C$ above depend on the
coefficients, which depend only on
the bound of $\hat Y$ and is independent of $n$, thanks to (\ref
{XYproperty}). Thus, letting $n\to\infty$ we have
\begin{eqnarray*}
\mathbb{E} \biggl\{\exp \biggl(\varepsilon \int_0^{\tau }
\llvert \hat Z_t\rrvert ^2\,dt \biggr) +\bigl\llvert
M^{\Lambda (\hat Y) \hat Z}_{\tau }\bigr\rrvert ^{1+\varepsilon
} \biggr\} \le C<
\infty.
\end{eqnarray*}
On the other hand, since $\cX$ satisfies (\ref{X}) on $[0, \tau )$,
we see that the estimate above implies that $\cX_\tau >0$, a.s. Thus,
$\tau = T$ a.s.
In other words, $(\hat Y, \hat Z)$ satisfies (\ref{hatYBSDE}) over
$[0, T]$, and (\ref{hatYproperty}) holds. The proof is now complete.
\end{pf}

%
\begin{rem}
We should point out that Theorems \ref{teo-linear-wellposed} only
indicates an {a priori} relationship between the characteristic
BSDE and the ``derivative'' of the decoupling field, whenever
exists, via the variational FBSDE (\ref{linearFBSDE}). The boundedness
requirement (\ref{hatYproperty}), or equivalently, the ``regularity''
of the decoupling field, is crucial for the
solution scheme to be effective (recall the inductive procedure in
Theorem \ref{teo-decoupling}). The actual
construction of the decoupling field, however, depends on the
well-posedness of the dominating ODEs
to be analyzed in details in next section, which is motivated by but
independent of the results in this section.
In fact, only a localized version (in small time duration) of Theorem
\ref{teo-linear-wellposed}(ii) will be used in the proof of Theorem
\ref{teo-smallT1}(iii) below.
\end{rem}

We conclude this section by presenting a result regarding the
uniqueness of the solutions to
FBSDE (\ref{linearFBSDE}) and its characteristic BSDE (\ref{hatYBSDE}),
which might be of independent interest. We should note that this result
will not be used in our future discussion, but
its arguments will be useful whenever a linearized FBSDE is encountered
(e.g., the proof of Theorem \ref{teo-comparison} below).
To this end,
we make use of an additional condition on $(\hat Y, \hat Z)$ that strengthen
the estimate (\ref{BMO}):
%
\begin{eqnarray}
\label{hatZproperty} \mathbb{E} \Bigl\{\sup_{0\le t\le T} \bigl\llvert
M^{\Lambda (\hat Y)\hat
Z}_t\bigr\rrvert ^{2+\varepsilon } \Bigr\}<\infty\qquad\mbox{for some }\varepsilon >0.
\end{eqnarray}

\begin{teo}
\label{teo-linearuniq} Let Assumption \ref{assum-linear} hold.
Then the BSDE (\ref{hatYBSDE}) has a solution
$(\hat Y, \hat Z)$
satisfying (\ref{hatYproperty}) and (\ref{hatZproperty}) if and only if
the FBSDE (\ref{linearFBSDE}) has a solution
$(\cX, \cY, \cZ)\in\widehat{\mathbb{L}}_2$ satisfying (\ref{XYproperty}).

Moreover, in such a case the uniqueness holds for solutions to BSDE
(\ref{hatYBSDE}) satisfying (\ref{hatYproperty}) and (\ref
{hatZproperty}) with $\varepsilon =0$ and for solutions to
FBSDE (\ref{linearFBSDE}) in $\mathbb{L}^2$ satisfying~(\ref{XYproperty}).
\end{teo}

\begin{pf}
We proceed in three steps. To make the presentation more
precise, we denote:
\begin{itemize}
\item $\mathcal{A}_0 \stackrel{\triangle}{=}\{\mbox{all solutions
$(\cX, \cY, \cZ)\in \mathbb{L}^2$ to FBSDE (\ref{linearFBSDE})
satisfying (\ref{XYproperty})}\}$;

\item$ \mathcal{A}\stackrel{\triangle}{=}\mathcal{A}_0 \cap
\widehat{\mathbb{L}}^2$;

\item $ \mathcal{B}_0 \stackrel{\triangle}{=}\{$all solutions
$(\hat Y, \hat Z)$ to BSDE (\ref{hatYBSDE}) satisfying (\ref
{hatYproperty}) and
(\ref{hatZproperty}) with $\varepsilon =0\}$; and

\item$\mathcal{B}\stackrel{\triangle}{=}\{$all solutions
$(\hat Y, \hat Z)$ in $\mathcal{B}_0$ satisfying (\ref
{hatZproperty})$\}$.
\end{itemize}

\begin{longlist}[\textit{Step}~2.]
\item[\textit{Step} 1.] We first prove the equivalence of the existence of
desired solutions in~$\mathcal{A}$ and $\mathcal{B}$. First, assume
there exists
$(\hat Y, \hat Z)\in\mathcal{B}$. Then\vspace*{1pt} by Theorem \ref{teo-linear-wellposed},
the \mbox{FBSDE (\ref{linearFBSDE})} has a solution $(\cX, \cY, \cZ)\in
\widehat{\mathbb{L}}_1$. Furthermore, using condition (\ref
{hatZproperty}) we can
actually improve the estimates (\ref{X1eps}) and (\ref{YZ1eps}) to
$\mathbb{L}^{2+\varepsilon /2}$, and thus \mbox{$(\cX, \cY, \cZ)\in
\mathcal{A}$.}

Conversely, there exists $(\cX, \cY, \cZ)\in\mathcal{A}\subseteq
\widehat{\mathbb{L}}_1$, then by Theorem \ref{teo-linear-wellposed}(ii), the $(\hat Y,
\hat Z)$ defined by (\ref{hatYZ-linear}) satisfy (\ref{hatYBSDE}) and
(\ref{hatYproperty}), and $\cX$ satisfy (\ref{X}). Thus,
\begin{eqnarray*}
M^{\Lambda (\hat Y) \hat Z}_t &=& \cX_t \bigl[M_t^{I(\hat Y)}
\bigr]^{-1}\exp \biggl\{-\int_0^t
\bigl[H_s(\hat Y_s, \hat Z_s)-I_s(
\hat Y_s)\Lambda _s(\hat Y_s) \hat
Z_s\bigr]\,ds \biggr\}.
\end{eqnarray*}
If\vspace*{1pt} $\mathbb{E}\{\sup_{0\le t\le T}\llvert  \cX_t\rrvert  ^{p}\}<\infty$ for some
$p>2$, by
estimates similar to (\ref{X1eps}) we obtain~(\ref{hatZproperty}).

\item[\textit{Step} 2.] We next turn to the uniqueness. We claim that:
%
\begin{equation}
\label{linearuniq-claim}
\quad \begin{tabular}{p{290pt}}
For any $(\cX, \cY, \cZ)\in \mathcal{A}_0$ and $(\hat Y, \hat Z)\in\mathcal{B}_0$,
if either $(\cX, \cY, \cZ)\in \mathcal{A}$ or $(\hat Y, \hat Z)\in\mathcal{B}$,
then relation (\ref{hatYZ-linear}) and equivalently (\ref{YZ-linear}) must
hold.
\end{tabular} %
\end{equation}
Now fix an $(\cX^0, \cY^0, \cZ^0)\in \mathcal{A}$ and $(\hat Y^0,
\hat Z^0)\in\mathcal{B}$ which satisfy (\ref{hatYZ-linear}) and
(\ref{YZ-linear}). For any $(\hat Y, \hat Z)\in\mathcal{B}_0$, apply
(\ref{linearuniq-claim}) on $(\cX^0, \cY^0, \cZ^0)$ and $(\hat Y,
\hat Z)$, we see that they satisfy (\ref{hatYZ-linear}), and thus
$(\hat Y, \hat Z)=(\hat Y^0, \hat Z^0)$.
On the other hand, for any $(\cX, \cY, \cZ)\in \mathcal{A}_0$,
apply (\ref{linearuniq-claim}) on $(\cX, \cY, \cZ)$ and $(\hat Y^0,
\hat Z^0)$. By\vspace*{1pt}
(\ref{YZ-linear}), we see that $\cX$ must satisfy (\ref{X}) with
$(\hat Y^0, \hat Z^0)$ in the right-hand\vspace*{1pt} side, and thus $\cX= \cX^0$.
Moreover, it follows from
(\ref{YZ-linear}) that $(\cY, \cZ)= (\cY^0, \cZ^0)$.

\item[\textit{Step} 3.] We now prove claim (\ref{linearuniq-claim}). Given $(\cX, \cY, \cZ)\in\mathcal{A}_0$ and $(\hat Y, \hat Z)\in\mathcal
{B}_0$, denote
%
\begin{eqnarray}\label{tildeYZ}
\delta Y_t &\stackrel{\triangle} {=}&\cY_t -
\hat Y_t \cX_t,
\nonumber\\[-8pt]\\[-8pt]\nonumber
\delta Z_t &\stackrel{\triangle} {=}&\cZ_t - \bigl[\cX_t \hat Z_t+
\hat Y_t(\sigma _1\cX_t+\sigma _2
\cY_t+\sigma _3\cZ_t) \bigr].
\end{eqnarray}
Applying It\^o's formula to $\delta Y_t$, we have
%
\begin{eqnarray}
\label{dDY0} d(\delta Y_t) &=& - \bigl[f_1
\cX_t + f_2 \cY_t + f_3
\cZ_t - \cX _t\bigl[F_t(\hat Y_t)
+ G_t(\hat Y_t) \hat Z_t + \Lambda
_t(\hat Y_t) \llvert \hat Z_t\rrvert
^2\bigr]
\nonumber
\\
&&\hspace*{40pt}{}+ \hat Y_t(b_1\cX_t+b_2
\cY_t+b_3\cZ_t) +\hat Z_t (\sigma
_1\cX _t+\sigma _2\cY_t+\sigma
_3\cZ_t) \bigr]\,dt\nonumber\hspace*{-20pt}
\\
&&{}  +\delta Z_t \,dB_t
\nonumber\\[-8pt]\\[-8pt]\nonumber
&=& - \bigl[\cX_t \bigl[f_1 - F_t(\hat
Y_t) - G_t(\hat Y_t) \hat Z_t -
\Lambda _t(\hat Y_t) \llvert \hat Z_t\rrvert
^2 + b_1 \hat Y_t +\sigma _1
\hat Z_t\bigr]
\\
&&\hspace*{51pt}{} + \cY_t [f_2 + b_2 \hat Y_t +
\sigma _2 \hat Z_t ] + \cZ _t
[f_3 + b_3\hat Y_t + \sigma _3
\hat Z_t] \bigr]\,dt\nonumber
\\
&&{} +\delta Z_t \,dB_t. \nonumber
\end{eqnarray}
By (\ref{tildeYZ}), one can easily
check that
\begin{eqnarray*}
\cY &=& \delta Y + \hat Y \cX,
\\
\cZ &=& {\delta Z + \cX\hat Z +
\hat Y(\sigma _1 \cX+ \sigma _2 \cY)\over1-\sigma _3 \hat
Y}
\\
&=&
{\delta Z+ \sigma _2 \hat Y \delta Y
+ \cX[\hat Z+ (\sigma _1+\sigma _2 \hat Y) \hat Y] \over1-\sigma
_3 \hat Y}.
\end{eqnarray*}
Plugging these into (\ref{dDY0}), we obtain
\begin{eqnarray*}
d(\delta Y_t) = - [\alpha _t \cX_t +
\beta_t\delta Y_t +\gamma _t \delta
Z_t ] \,dt +\delta Z_t \,dB_t,
\end{eqnarray*}
where
%
\begin{eqnarray}
\label{betagamma} \gamma _t &\stackrel{\triangle} {=}&
{f_3 + b_3\hat Y_t
+ \sigma _3 \hat Z_t\over1-\sigma _3\hat Y_t} = {f_3 + b_3\hat
Y_t\over1-\sigma _3\hat Y_t} + \Lambda _t(\hat
Y_t) \hat Z_t;
\nonumber\\[-8pt]\\[-8pt]\nonumber
\beta_t &\stackrel{\triangle} {=}& f_2 + b_2
\hat Y_t +\sigma _2 \hat Z_t +
{\sigma _2\hat Y_t
[f_3 + b_3\hat Y_t + \sigma _3 \hat Z_t]\over1-\sigma _3\hat Y_t}
\end{eqnarray}
and
\begin{eqnarray*}
\alpha _t &\stackrel{\triangle} {=}& f_1 -
F_t(\hat Y_t) - G_t(\hat Y_t)
\hat Z_t - \Lambda _t(\hat Y_t) \llvert \hat
Z_t\rrvert ^2 + b_1 \hat Y_t +
\sigma _1 \hat Z_t
\\
&&{}+ [f_2 + b_2 \hat Y_t + \sigma
_2 \hat Z_t] \hat Y_t + [f_3 +
b_3\hat Y_t + \sigma _3 \hat Z_t]
{\hat Z_t + (\sigma _1+\sigma _2
\hat Y_t) \hat Y_t \over1-\sigma _3 \hat Y_t}
\\
&=& 0,
\end{eqnarray*}
thanks to (\ref{F}). Denote
%
\begin{eqnarray}\label{Gamma}
\Gamma _t &\stackrel{\triangle} {=}& M^{\gamma }_t
\exp \biggl(\int_0^t \beta_s \,ds
\biggr)\nonumber
\\
&=& M^{\Lambda (\hat Y) \hat Z}_t M^{(f_3 + b_3\hat Y)/(1-\sigma _3\hat Y)}_t
\\
&&{}\times \exp
\biggl(\int_0^t \biggl[\beta_s -
{f_3 + b_3\hat Y_s\over1-\sigma
_3\hat Y_s}\Lambda _s(\hat Y_s)\hat
Z_s\biggr]\,ds \biggr).\nonumber
\end{eqnarray}
Then by applying It\^{o}'s formula, one obtains immediately
%
\begin{eqnarray}
\label{dDY} d(\Gamma _t \delta Y_t) &=& \Gamma
_t[\gamma _t\delta Y_t+\delta
Z_t]\,dB_t.
\end{eqnarray}
We claim that
%
\begin{eqnarray}\label{mgest3}
&&  \mathbb{E} \biggl\{ \biggl(\int_0^T
\llvert \Gamma _t\rrvert ^2[\gamma _t\delta
Y_t+\delta Z_t]^2\,dt \biggr)^{1/2}
\biggr\}
\nonumber\\[-8pt]\\[-8pt]\nonumber
&&\qquad \le  \mathbb{E} \biggl\{\sup_{0\le t\le T}\llvert \Gamma
_t\rrvert \biggl(\int_0^T [\gamma
_t\delta Y_t+\delta Z_t]^2\,dt
\biggr)^{1/2} \biggr\}<\infty,
\end{eqnarray}
so that $\int_0^\cdot \Gamma _s [\gamma _t\delta Y_t+\delta Z_s]
\,dB_s$ is a true martingale. Since $\delta Y_T=0$ and $\Gamma _0=1$,
it follows from (\ref{dDY})
that $\delta Y=0$, and hence $\delta Z=0$. Then (\ref{tildeYZ})
leads to (\ref{hatYZ-linear}) immediately.

It remains to prove (\ref{mgest3}). Note that
\begin{eqnarray*}
\llvert \gamma _t\rrvert &\le& C\bigl[1+\llvert \hat Z_t
\rrvert \bigr],\qquad \llvert \delta Y_t\rrvert \le C\bigl[\llvert
\cX_t\rrvert + \llvert \cY_t\rrvert \bigr],
\\
\llvert
\delta Z_t\rrvert &\le& C\bigl[\llvert \cX_t\rrvert +\llvert
\cY_t\rrvert +\llvert \cZ_t\rrvert + \llvert \cX
_t\rrvert \llvert \hat Z_t\rrvert \bigr].
\end{eqnarray*}
Then
%
\begin{eqnarray}
\label{mgest1} &&\int_0^T [\gamma
_t\delta Y_t+\delta Z_t]^2\,dt\nonumber
\\
&&\qquad \le C \biggl[1+ \sup_{0\le t\le T}\bigl[\llvert \cX_t
\rrvert ^2 + \llvert \cY_t\rrvert ^2\bigr]
\\
&&\hspace*{45pt}{} +
\int_0^T \bigl[\llvert \cZ_t\rrvert
^2+ \llvert \hat Z_t\rrvert ^2\bigr]\,dt + \sup
_{0\le t\le T}\llvert \cX_t\rrvert ^2 \int
_0^T \llvert \hat Z_t\rrvert
^2\,dt \biggr].\nonumber
\end{eqnarray}
Since $(\hat Y, \hat Z)$ satisfies (\ref{hatYproperty}), by (\ref
{BMO}) we have
%
\begin{eqnarray}
\label{mgest2} \mathbb{E} \biggl\{ \biggl(\int_0^T
\llvert \hat Z_t\rrvert ^2\,dt \biggr)^p \biggr
\}<\infty\qquad\mbox{for any } p \ge1.
\end{eqnarray}
We now verify (\ref{mgest3}) in the two cases:\vspace*{1pt}
\end{longlist}

\begin{longlist}[\textit{Case} 1.]
\item[\textit{Case} 1.] $(\hat Y, \hat Z)\in\mathcal{B}$, namely (\ref
{hatZproperty}) holds with some $\varepsilon > 0$. Following the
arguments for (\ref{X1eps}), we have
%
\begin{eqnarray}
\label{Gammaeps} \mathbb{E}\Bigl\{\sup_{t\in[0,T]}\llvert \Gamma
_t\rrvert ^{2+{\varepsilon /2}}\Bigr\} <\infty.
\end{eqnarray}
Then for any $(\cX, \cY, \cZ) \in\mathcal{A}_0$, plugging (\ref
{mgest2}) and (\ref{Gammaeps}) into (\ref{mgest1}) we have (\ref
{mgest3}) immediately.

\item[\textit{Case} 2.] $(\cX, \cY, \cZ) \in\mathcal{A}$, namely $(\cX, \cY, \cZ)\in\mathbb{L}^{2+\varepsilon }$ for some $\varepsilon >0$.
Note that
\begin{eqnarray*}
{1\over2+\varepsilon } + {3+2\varepsilon \over6 + 3\varepsilon } + {\varepsilon \over6+3\varepsilon } =
1\quad\mbox{and}\quad {6+3\varepsilon \over3+2\varepsilon } < 2.
\end{eqnarray*}
Since (\ref{hatZproperty}) holds with $\varepsilon =0$, following the
arguments for (\ref{X1eps}) we
have
$\mathbb{E}\{\sup_{t\in[0,T]}\llvert  \Gamma _t\rrvert  ^{(6+3\varepsilon)/(3+2\varepsilon)}\}<\infty$.
This implies that
\begin{eqnarray*}
&&\mathbb{E} \biggl\{\sup_{0\le t\le T} \llvert \Gamma _t
\rrvert \sup_{0\le t\le
T}\llvert \cX_t\rrvert \int
_0^T \llvert \hat Z_t\rrvert
^2\,dt \biggr\}
\\
&&\qquad \le\Bigl(\mathbb{E} \Bigl\{\sup_{0\le t\le T} \llvert \Gamma
_t\rrvert ^{(6+3\varepsilon )/(3+2\varepsilon) } \Bigr\} \Bigr)^{(3+2\varepsilon )/(6 + 3\varepsilon) }
\\
&&\quad\qquad{} \times \Bigl(
\mathbb{E} \Bigl\{ \sup_{0\le t\le T} \llvert \cX_t\rrvert
^{2+\varepsilon } \Bigr\} \Bigr)^{1/(2 + \varepsilon) } \biggl(\mathbb{E} \biggl\{ \biggl(
\int_0^T \llvert \hat Z_t\rrvert
^2\,dt \biggr)^{(6+3\varepsilon )/\varepsilon } \biggr\} \biggr)^{\varepsilon /(6 + 3\varepsilon )}
\\
&&\qquad <\infty.
\end{eqnarray*}
Then one can easily prove (\ref{mgest3}) again.\quad\qed
\end{longlist}\noqed
\end{pf}

\section{Well-posedness of the dominating equations}\label{sec5}
\label{sect-ODE}

We note that Theorems \ref{teo-linear-wellposed} and \ref{teo-linearuniq} only established the relations of the well-posedness
between the characteristic
BSDEs and the original FBSDE, it does not provide the well-posedness
result for either one of them.
In this section, we take a closer look at the dominating ODEs
(\ref{ODE}). Since the existence of bounded solutions $\overline\by$
and $\underline\by$ to the dominating ODEs will be
essential in constructing the desired regular decoupling field, 
which will eventually lead to well-posedness of the FBSDE (\ref
{FBSDE}), the results in this section will be the blueprint of a
\textit{user's guide} in the end.

We begin with a special form of comparison theorem among the solutions
to ODEs.
Consider the following ``backward ODEs'' on $[0,T]$:
%
\begin{eqnarray}
\label{ODE0} \by^0_t = h^0 + \int
_t^T F^0\bigl(s,
\by^0_s\bigr) \,ds
\end{eqnarray}
and
%
\begin{eqnarray}\label{ODEoly}
\by^1_t &=& h^1 -
C^1+ \int_t^T\bigl[F^1
\bigl(s,\by^1_s\bigr)+c^1_s
\bigr]\,ds,
\nonumber\\[-8pt]\\[-8pt]\nonumber
\by^2_t &=& h^2 + C^2+
\int_t^T\bigl[F^2\bigl(s,
\by^2_s\bigr)-c^2_s\bigr]\,ds,
\end{eqnarray}
where $F^0, F^1, F^2\dvtx [0,T]\times\mathbb{R}\mathop{\longrightarrow}\mathbb{R}$ are (deterministic) measurable functions.
The following simple lemma will be useful in our discussion. Its proof
is rather elementary and we defer it to the \hyperref[app]{Appendix}.

\begin{lem}
\label{lem-ODEcomp}
Assume that:
\begin{longlist}[(iii)]
\item[(i)] $h^1 \le h^0\le h^2$, and $F^1 \le F^0 \le F^2$.

\item[(ii)] Both ODEs in (\ref{ODEoly}) admit bounded solutions $\by^1$ and
$\by^2$ on $[0,T]$.

\item[(iii)] For\vspace*{1pt} each $t\in[0,T]$, the functions $y\mapsto F^i(t,y)$,
$i=0,1,2$, are uniformly Lipschitz continuous for
$y\in[\by^1_t, \by^2_t]$, with a common Lipschitz constant $L$.

\item[(iv)] $
C^i \ge\int_t^T e^{-\int_s^T \alpha _r \,dr} c^i_s\,ds$, for all $t\in
[0,T]$ and all $\alpha $ satisfying $\llvert  \alpha \rrvert  \le L$.
\end{longlist}
Then (\ref{ODE0}) has a unique solution $\by^0$ satisfying
$\by^1 \le\by^0 \le\by^2$.
\end{lem}

\begin{rem}
\label{rem-ODEcomp}
A typical sufficient condition for the above (iv) is: $C^i \ge
\int_0^T e^{L(T-t)} (c^i_t)^+\,dt$. In particular, this is satisfied if
$C^i=0$ and $c^i \le0$.
\end{rem}

\subsection{Linear FBSDE with constant coefficients}\label{sec5.1}

We first investigate the linear FBSDE (\ref{linearFBSDE}) where
all the coefficients are constants.
We shall show that in such a case some ``sharp''
(sufficient and necessary) conditions regarding well-posedness can be
obtained. These
results, to our best knowledge, are novel in the literature; and at the
same time, they more
or less set the ``limits'' for the solvability of general FBSDE (\ref{FBSDE}).

We carry out our analysis in two cases.

\begin{longlist}[\textit{Case} 1.]
\item[\textit{Case} 1: $\sigma _3=0$.]
In this case, $\overline h = \underline h = h$, $\overline F(t,y)=
\underline F(t,y) = F(y)$, and
two ODEs in (\ref{ODE}) become the same:
%
\begin{eqnarray}
\label{ODE1} \by_t = h + \int_t^T
F(\by_s) \,ds,
\end{eqnarray}
where
%
\begin{equation}
\label{F1} \qquad F(y) = f_1 + [f_2+b_1+
\sigma _1f_3]y+[b_2+f_3\sigma
_2+b_3\sigma _1]y^2+ \sigma
_2b_3y^3.
\end{equation}
%
We have the following theorem.
\end{longlist}

%
\begin{teo}
\label{teo-constant-F1} Assume that
in the linear FBSDE (\ref{linearFBSDE}) all coefficients are
constants, and
$\sigma _3=0$. Then the corresponding dominating ODE (\ref{ODE1})
with $F$ defined by (\ref{F1}) has a bounded solution for arbitrary
$T$ if and only if one of the following three cases hold
true:
\begin{longlist}[(iii)]
\item[(i)] $F(h)\ge0$ and $F$ has a zero point in $[h, \infty)$.

\item[(ii)] $F(h)\le0$ and $F$ has a zero point in $(-\infty, h]$.

\item[(iii)] $\sigma _2b_3=0$ and $b_2+f_3\sigma _2+b_3\sigma _1=0$.
\end{longlist}
\end{teo}

\begin{pf}
We first prove the sufficiency part. In case (i), there
exists $\lambda \ge h$ such that $F(\lambda )=0$. Note that $F$ is
locally Lipschitz continuous in $y$ and
\begin{eqnarray*}
h = h + \int_t^T \bigl[F(h) - F(h)\bigr]\,ds,
\qquad \lambda = \lambda + \int_t^T F(\lambda )
\,ds.
\end{eqnarray*}
Then it follows from Lemma \ref{lem-ODEcomp} and in particular Remark
\ref{rem-ODEcomp} that $\by_t\in[h,\lambda ]$, $t\in[0,T]$.
Similarly, in case (ii), one has $\by_t\in[\lambda, h]$, for some
$\lambda \le h$
such that $F(\lambda )=0$. Finally,
in case (iii) the
ODE (\ref{ODE1}) becomes linear:
%
\begin{eqnarray}
\by_t = h + \int_t^T
\bigl[f_1 + (f_2+b_1+\sigma
_1f_3)\by_s\bigr] \,ds.
\end{eqnarray}
Thus, it is obviously bounded.

The proof of necessity is elementary but lengthy,
we postpone it to the \hyperref[app]{Appendix}.
\end{pf}

When the terminal time $T$ is fixed, we have the following slightly
weaker sufficient conditions:

\begin{teo}
\label{teo-constant-F2} For any given $T>0$, the ODE (\ref{ODE1})
with $F$ given in (\ref{F1}) has a
bounded solution on $[0, T]$ if one of the following three cases hold true:
\begin{longlist}[(iii)]
\item[(i)] $\sigma _2 b_3 < 0$ or $F(h)=0$.

\item[(ii)] $F(h)>0$, and there exists a constant $\varepsilon =\varepsilon
(T)>0$ small enough, such
that
\[
\sigma _2 b_3 \le\varepsilon\quad\mbox{and}\quad
b_2+ f_3\sigma _2 + b_3\sigma _1\le\varepsilon .
\]

\item[(iii)] $F(h)< 0$, and there exists a constant $\varepsilon
=\varepsilon (T)>0$ small
enough such that
\[
\sigma _2 b_3 \le\varepsilon\quad\mbox{and}\quad
b_2 + f_3\sigma _2 + b_3\sigma _1\ge-\varepsilon.
\]
\end{longlist}
\end{teo}

\begin{pf}
(i) In this case clearly, the result follows from either
(i) or (ii) of Theorem \ref{teo-constant-F1}.

(ii) In this case we have, for some small constant $\varepsilon >0$
which will be specified later and for some constants $C_1, C_0$
independent of $\varepsilon $,
%
\begin{equation}
\label{Fupbound}\qquad F(y) \le\varepsilon y^3 + \varepsilon y^2
+ C_1 y + C \le2 \varepsilon y^3 + C_1y +
C_0\qquad\mbox{for all } y\ge0.
\end{equation}
We first solve
\begin{eqnarray*}
\tilde\by_t = h^+ + \int_t^T [
C_1\tilde\by_s + C_0 + 1]\,ds
\end{eqnarray*}
and obtain
%
\begin{eqnarray}
\label{C2}
\tilde\by_t &=& e^{C_1(T-t)}h^+ +
{C_0+1\over C_1}\bigl[e^{C_1(T-t)} -1\bigr]
\nonumber\\[-8pt]\\[-8pt]\nonumber
&\le& C_2 :=
e^{C_1T}h^+ + {C_0+1\over C_1}\bigl[e^{C_1T} -1\bigr].
\end{eqnarray}
Set $\varepsilon \stackrel{\triangle}{=}{1\over2C_2^3}$ so that
$2\varepsilon \tilde\by_t^3 \le1$. Note that
\begin{eqnarray*}
\tilde\by_t = h^+ + \int_t^T
\bigl[2\varepsilon \tilde\by_s^3 + C_1\tilde
\by_s + C_0 + \bigl(1-2\varepsilon \tilde
\by_s^3\bigr)\bigr]\,ds.
\end{eqnarray*}
By (\ref{Fupbound}), applying Lemma \ref{lem-ODEcomp} and in
particular Remark \ref{rem-ODEcomp} we see that ODE
(\ref{ODE1}) has a solution $\by\in[h, \tilde\by]\subset[h,
C_2]$.

(iii) can be proved similarly.
\end{pf}

\begin{longlist}[\textit{Case} 2.]
\item[\textit{Case} 2: $\sigma _3 \neq0$.]
In this case, we still have $\overline h = \underline h = h$,
$\overline F(t, y) = \underline F(t,y) = F(y)$, where the deterministic
function $F$ in (\ref{F1}) can be rewritten as
%
\begin{eqnarray}
\label{F2} F(y) = {\alpha _0 \over{1/\sigma _3}- y} + \alpha _1+\alpha
_2y + \biggl[b_2-{b_3\sigma _2\over\sigma _3}
\biggr]y^2,
\end{eqnarray}
for some constants $\alpha _0, \alpha _1, \alpha _2$. In this case,
the two ODEs in (\ref{ODE}) also become the same one (\ref{ODE1})
and, in light of (\ref{hatYproperty}), we want to find its solution
satisfying that
%
\begin{eqnarray}
\label{yproperty} \mbox{both $\by$ and $(1-\sigma _3
\by)^{-1}$ are bounded.}
\end{eqnarray}

\begin{rem}
\label{sigma3h}
We note that (\ref{yproperty}) amounts to saying that $\sigma _3
h\neq1$ since \mbox{$\by_T=h$}. In fact, if $\sigma _3 h=1$, there are
counter examples in both existence and uniqueness of the linear FBSDE
(\ref{linearFBSDE}) (cf., e.g., \cite{mybk}).
\end{rem}

We now have
the following theorem.

\begin{teo}
\label{teo-constant-F3}
Assume the FBSDE is the linear one (\ref{linearFBSDE}) and all the
coefficients are constants. Assume also that $\sigma _3\neq0$ and
$h\sigma _3\neq1$.
Then the ODE (\ref{ODE1}) has a solution satisfying (\ref{yproperty}) for
arbitrary $T$ if and only if one of the following four cases holds:
\begin{longlist}[(iii)]
\item[(i)] $h < {1\over\sigma _3}$, $F(h) \le0$, and either $F$ has a zero
point in $(-\infty, h]$ or \mbox{$b_2-{b_3\sigma _2\over\sigma _3}=0$}.

\item[(ii)] $h > {1\over\sigma _3}$, $F(h) \ge0$,\vspace*{1pt} and either $F$ has a zero
point in $[h, \infty)$ or $b_2-{b_3\sigma _2\over\sigma _3}=0$.

\item[(iii)] $h < {1\over\sigma _3}$, $F(h) \ge0$,\vspace*{1pt} and $F$ has a zero point
in $[h, {1\over\sigma _3})$.

\item[(iv)] $h >
{1\over\sigma _3}$, $F(h) \le0$, and $F$ has a zero point in
$({1\over\sigma _3}, h]$.
\end{longlist}
\end{teo}

\begin{pf}
We prove the sufficiency here and again postpone the
necessary part to the \hyperref[app]{Appendix}.
\begin{longlist}[(iii)]
\item[(i)] If $F(\lambda )=0$, for some $\lambda \in(-\infty, h]$, then as
in Theorem \ref{teo-constant-F1} we see that ODE (\ref{ODE1}) has a\vspace*{1pt}
solution $\by\in[\lambda, h]$. Thus (\ref{yproperty}) holds. We
now assume instead that $b_2-{b_3\sigma _2\over\sigma _3}=0$. Then\vspace*{1pt}
from (\ref{F2}), we see that
$F(y)=\alpha _0 ({1\over\sigma _3}- y)^{-1}+\alpha _1+\alpha _2y$. Consider
\begin{eqnarray*}
\tilde\by_t = h + \int_t^T
\biggl[-\llvert \alpha _0\rrvert \biggl({1\over\sigma _3}- h
\biggr)^{-1}+\alpha _1+\alpha _2\tilde
\by_s\biggr] \,ds.
\end{eqnarray*}
Since $F(h) \le0$, clearly the above SDE has a bounded solution
$\tilde\by\le h$. Applying Lemma \ref{lem-ODEcomp}, one can easily
see that (\ref{ODE1}) has a solution $\by\in[\tilde\by, h]$. Thus,
(\ref{yproperty}) holds.

\item[(iii)] Let $\lambda \in[h,{1\over\sigma _3})$ be such that
$F(\lambda )=0$. Note that $\by^1_t \stackrel{\triangle}{=}h$ and
$\by^2_t \stackrel{\triangle}{=}\lambda $ are (constant) solutions
of the following ODEs, respectively:
%
\begin{eqnarray*}
\by^1_t=h+\int_t^T
\bigl[F\bigl(\by^1_s\bigr)-F(h)\bigr]\,ds,\qquad
\by^2_t=\lambda +\int_t^T
F\bigl(\by^2_s\bigr)\,ds.
\end{eqnarray*}
Comparing these two equations with (\ref{ODE1}) and applying Lemma
\ref{lem-ODEcomp}, we have $h\le\by_t\le\lambda $, for any $t\in
[0,T]$. This implies (\ref{yproperty}) immediately.

(ii) and (iv) can be proved similarly as (i) and (iii), respectively.\quad\qed
\end{longlist}\noqed
\end{pf}

When $T$ is fixed, we may also have some slightly weaker sufficient
conditions. However, these conditions are more involved, so we omit
them here and will discuss directly for the general case in next
subsection; see Theorems \ref{teo-nonlinear2} and \ref{teo-nonlinear3} below.
\end{longlist}

\subsection{The nonlinear case}\label{sec5.2}

Again we consider the case that $\sigma _3 = 0$ first.
\begin{longlist}[\textit{Case} 2.]
\item[\textit{Case} 1: $\sigma = \sigma (t, x, y)$.]
We recall that in this case $F$ takes the form (\ref{F1}), where
$b_i$, $\sigma _i$, $f_i$, $i=1,2,3$, are bounded, adapted processes
defined by (\ref{tildephi}), and thus $F$ is also random and may depend
on $t$. Now recall the definition of the functions $\overline F$ and
$\underline
F$ in (\ref{Fbar}). Again, by a slight abuse of notation we replace
$\Theta ^j$, $j=1,2$ in (\ref{tildephi}) by $\theta _j$, $j=1,2$,
and still
denote them by $b_i$, $\sigma _i$, $f_i$, $i=1,2,3$. In what follows, all
assumptions involving coefficients in (\ref{F1}) will be in the
sense that they hold uniformly for all $\theta _j$, $j=1,2$. In analogy
to Theorem \ref{teo-constant-F2}, we have the following result.
\end{longlist}

\begin{teo}
\label{teo-nonlinear1} Assume Assumption
\ref{assum-Lipschitz} holds and $\sigma = \sigma (t, x, y)$. Then,
for any $T>0$, the ODEs
(\ref{ODE}) have bounded solutions $\overline\by$ and $\underline
\by$ on $[0,T]$
if one of the following three cases holds true:
\begin{longlist}[(iii)]
\item[(i)] There exists a constant $\varepsilon >0$ such that
%
\begin{eqnarray}
\label{nonlinear-small0} \sigma _2 b_3 \le-\varepsilon \llvert
b_2 + f_3\sigma _2 + b_3\sigma
_1\rrvert.
\end{eqnarray}

\item[(ii)] There exists a constant
$\lambda \le\underline h$, and a constant $\varepsilon >0$ small
enough such
that
%
\begin{eqnarray}
\label{nonlinear-small1} \underline F(t, \lambda ) \ge0,\qquad \sigma _2
b_3 \le\varepsilon\quad\mbox{and}\quad b_2 + f_3
\sigma _2 + b_3\sigma _1\le\varepsilon.
\end{eqnarray}

\item[(iii)] There exists a constant
$\lambda \ge\overline h$, and a constant $\varepsilon >0$ small
enough such
that
%
\begin{eqnarray}
\label{nonlinear-small2} \overline F(t, \lambda )\le 0,\qquad \sigma _2
b_3 \le\varepsilon\quad\mbox{and}\quad b_2 + f_3
\sigma _2 + b_3\sigma _1\ge-\varepsilon.
\end{eqnarray}
\end{longlist}
\end{teo}

\begin{pf}
(i) In this case, we have
\begin{eqnarray*}
\overline F(t, y) &\le& C[y + 1]\qquad\mbox{for all }y\ge{1\over
\varepsilon }\quad\mbox{and}
\\
\underline F(t, y) &\ge& C[y - 1]\qquad\mbox{for all }y\le-
{1\over\varepsilon }.
\end{eqnarray*}
Following the arguments in Theorem \ref{teo-constant-F2}(ii), one can
easily prove the result.

(ii) In this case, similar to (\ref{Fupbound}) we have
\begin{eqnarray}
\overline F(t,\lambda ) \ge\underline F(t,\lambda ) \ge0\quad\mbox{and}\quad
\underline F(t,y) \le\overline F(t,y) \le2\varepsilon y^3 +
C_1 y + C_0
\nonumber\\
\eqntext{\mbox{for all }y\ge0.}
\end{eqnarray}
Let $C_2$ be defined by (\ref{C2}) and set $\varepsilon:= {1\over2
C_2^3}$. Following the arguments in Theorem~\ref{teo-constant-F2}(ii),
we see that the ODEs in (\ref{ODE}) have bounded solutions
$\lambda \le\underline \by\le\overline \by\le C_2$.

(iii) can be proved similarly.
\end{pf}

\begin{longlist}
\item[\textit{Case} 2: $\sigma = \sigma (t,x,y,z)$.]
This case has been avoided in many of the existing literature,
especially when one uses the decoupling strategy.
A well-known sufficient condition for the existence is, roughly
speaking, that
$
\llvert  \sigma _3h\rrvert  <1$.
As we will see below, the condition we need is essentially $\sigma _3
h \neq1$. In particular, we shall discuss three different cases:
\begin{longlist}[(2-a)]
\item[(2-a)] $\llvert  \sigma _3h\rrvert  <1$;

\item[(2-b)] $\llvert  \sigma _3 h\rrvert   >1$ and both $\sigma _3$ and $h$ do not change sign;

\item[(2-c)] $\sigma _3 h <1$ and either $\sigma _3$ or $h$ does not change sign.
\end{longlist}
\end{longlist}

\begin{rem}
\label{linearcase}
We remark that, if all the coefficients are constants, the above three
cases (actually the latter two) cover all possible cases of $\sigma _3
h\neq1$. However, for general nonlinear FBSDEs with random
coefficients, we need them to hold uniformly in certain sense.
\end{rem}

To be more precise, let $T>0$ be given. We begin by fixing three
constants $c_1, c_2, c_3$ satisfying
%
\begin{eqnarray}
\label{c} c_1> 0,\qquad 0< c_2 < c_3,
\qquad c_1 c_3 <1.
\end{eqnarray}
The following result gives the answer to case (2-a).

\begin{teo}
\label{teo-nonlinear2} Assume that Assumption
\ref{assum-Lipschitz} and (\ref{c}) are in force. Assume also that
there exists a constant $\varepsilon = \varepsilon (T)>0$ small
enough such that
%
\begin{equation}
\label{case1} \qquad\llvert \sigma _3\rrvert \le c_1,\qquad
\llvert h\rrvert \le c_2\quad\mbox{and}\quad \overline F(t,c_3)
\le\varepsilon,\qquad \underline F(t,-c_3)\ge-\varepsilon.
\end{equation}
Then the ODEs in (\ref{ODE}) have solutions $\overline\by$ and
$\underline\by$ satisfying
\begin{eqnarray*}
&& -c_3 \le\underline\by\le\overline\by\le c_3\quad\mbox{and hence}\quad\mbox{both $\overline\by$ and $\underline\by$ satisfy (\ref{yproperty})}.
\end{eqnarray*}
\end{teo}

\begin{pf}
Note that $1- \sigma _3 y \ge1- c_1 c_3>0$ for $y\in
[-c_3, c_3]$, then $\overline F$ and $\underline F$ are uniformly Lipschitz
continuous in $y$ for $y\in[-c_3, c_3]$, and we denote by $L$ their
uniform Lipschitz constant. Clearly, $\tilde\by^1_t \stackrel
{\triangle}{=}- c_3$
and $\tilde\by^2_t \stackrel{\triangle}{=}c_3$ satisfy the
following ODEs:
\begin{eqnarray*}
\tilde\by^1_t &=& -c_2- (c_3 -
c_2) + \int_t^T \bigl[\underline F
\bigl(s, \tilde\by^1_s\bigr) - \underline F(s,
-c_3)\bigr]\,ds,
\\
\tilde\by^2_t &=& c_2 + (c_3 -
c_2) + \int_t^T \bigl[\overline F
\bigl(s, \tilde\by^2_s\bigr) - \overline F(s,
c_3)\bigr]\,ds.
\end{eqnarray*}
Now set $\varepsilon >0$ small enough such that
$
c_3 - c_2 > \int_0^T e^{L(T-t)} \varepsilon \,dt$.
Then it follows from Lemma \ref{lem-ODEcomp} and in particular Remark
\ref{rem-ODEcomp} we obtain the result.
\end{pf}

We next consider case (2-b).

\begin{teo}
\label{teo-nonlinear3}
Let Assumption \ref{assum-Lipschitz} and (\ref{c}) hold. Assume that
there exists a constant $\varepsilon >0$ small enough such that one of
the following four cases holds true:
%
\begin{eqnarray}
\sigma _3&\ge& c_1^{-1}, \qquad h
\ge c_2^{-1}\quad\mbox{and}
\nonumber\\[-8pt]\label{case21}  \\[-8pt]\nonumber
\underline F\bigl(t,c^{-1}_3 \bigr)&\ge&-\varepsilon,\qquad b_2 - {b_3\sigma _2\over\sigma _3} \le
\varepsilon;
\\
\sigma _3&\le&-c_1^{-1}, \qquad h
\ge c_2^{-1}\quad\mbox{and}
\nonumber\\[-8pt]\label{case22}  \\[-8pt]\nonumber
\underline F\bigl(t,c^{-1}_3
\bigr)&\ge& -\varepsilon,\qquad b_2 - {b_3\sigma _2\over\sigma _3} \le
\varepsilon;
\\
\sigma _3&\ge& c_1^{-1}, \qquad h
\le- c_2^{-1}\quad\mbox{and}
\nonumber\\[-8pt]\label{case23} \\[-8pt]\nonumber
\overline F\bigl(t,-c^{-1}_3
\bigr)&\le&\varepsilon,\qquad b_2 - {b_3\sigma
_2\over\sigma _3} \ge -
\varepsilon;
\\
\sigma _3&\le& -c_1^{-1}, \qquad h
\le- c_2^{-1}\quad\mbox{and}
\nonumber\\[-8pt]\label{case24} \\[-8pt]\nonumber
\overline F\bigl(t,-c^{-1}_3
\bigr)&\le&\varepsilon,\qquad b_2 - {b_3\sigma
_2\over\sigma _3} \ge -
\varepsilon.
\end{eqnarray}
Then the ODEs in (\ref{ODE}) have bounded solutions $\overline\by$
and $\underline\by$ such that they satisfy the corresponding property
of $h$ in the above conditions with $c_2$ being replaced by $c_3$. In
particular, both $\overline\by$ and $\underline\by$ satisfy (\ref
{yproperty}).
\end{teo}

\begin{pf}
We prove only the case (\ref{case21}). The other cases
can be proved similarly.

In this case, we have
\begin{eqnarray*}
\overline F(t,y) \le{C\over c_3^{-1} - c_1} + C_1 y + \varepsilon
y^2 = C_0 + C_1 y + \varepsilon
y^2\qquad\mbox{for all }y \ge c_3^{-1}.
\end{eqnarray*}
Let $\tilde\by$ denote the bounded solution to the following ODE:
\begin{eqnarray*}
\tilde\by_t = \overline h + \int_t^T
[C_1 \tilde\by_s + C_0 + 1]\,ds\quad\mbox{and}\quad
C_2:= \tilde\by_0=\sup_{0\le t\le T} \tilde
\by_t.
\end{eqnarray*}
Let $L$ denote the uniform Lipschitz constant of $\underline F$ and
$\overline F$ for $y\in[c_3^{-1}, C_2]$. Note that
$
\underline F(t, c_3^{-1}) \ge-\varepsilon $.
Now follow the arguments in Theorem \ref{teo-nonlinear2} for the lower
bound and those in Theorem \ref{teo-constant-F2}(ii) for the upper
bound, one can easily show that, for $\varepsilon $ sufficiently
small, the ODEs in (\ref{ODE}) have solutions $\overline\by$ and
$\underline\by$ such that $c_3^{-1}\le\underline\by\le\overline
\by\le C_2$.
\end{pf}

We finally present the result for case (2-c).

\begin{teo}
\label{teo-nonlinear4}
Let Assumption \ref{assum-Lipschitz} and (\ref{c}) hold. Assume there
exists a~constant $\varepsilon >0$ small enough such that one of the
following four cases holds true:
%
\begin{eqnarray}
\sigma _3 &\le& c_1, \qquad 0\le h\le
c_2\quad\mbox{and}
\nonumber\\[-8pt]\label{case31} \\[-8pt]\nonumber
\overline F(t,c_3)&\le&\varepsilon, \qquad f_1 \ge0;
\\
0&\le&\sigma _3\le c_1,\qquad h\le
c_2\quad\mbox{and}
\nonumber\\[-8pt]\label{case32} \\[-8pt]\nonumber
\overline F(t,c_3)&\le&\varepsilon,\qquad
b_2 - {b_3\sigma _2\over\sigma _3} \ge-\varepsilon;
\\
\sigma _3 &\ge&-c_1, \qquad 0\ge h
\ge-c_2\quad\mbox{and}
\nonumber\\[-8pt]\label{case33} \\[-8pt]\nonumber
\underline F(t,-c_3)&\ge& -\varepsilon,\qquad f_1 \le0;
\\
0&\ge&\sigma _3 \ge-c_1,\qquad h
\ge-c_2\quad\mbox{and}
\nonumber\\[-8pt]\label{case34} \\[-8pt]\nonumber
\underline F(t,-c_3)&\ge&-\varepsilon,
\qquad b_2 - {b_3\sigma _2\over\sigma
_3} \le\varepsilon.
\end{eqnarray}
Then the ODEs in (\ref{ODE}) have bounded solutions $\overline\by$
and $\underline\by$ such that they satisfy the corresponding property
of $h$ in the above conditions with $c_2$ being replaced by $c_3$. In
particular, both $\overline\by$ and $\underline\by$ satisfy (\ref
{yproperty}).
\end{teo}

\begin{pf}
If (\ref{case31}) holds, then
$\underline F(t, 0) \ge0$ and $F(t, c_3) \le\varepsilon $.
Following the arguments in Theorem \ref{teo-constant-F1} for the lower
bound and those in Theorem \ref{teo-nonlinear1} for the upper bound,
one can easily show that, for $\varepsilon $ sufficiently small, the
ODEs in (\ref{ODE}) have solutions $\overline\by$ and $\underline
\by$ such that $0\le\underline\by\le\overline\by\le c_3$.

If (\ref{case32}) holds, follow the arguments in Theorem \ref{teo-constant-F2}(ii) for the lower
bound and those in Theorem \ref{teo-nonlinear1} for the upper bound, one can easily show that, for
$\varepsilon $ sufficiently small, the ODEs in (\ref{ODE}) have
solutions $\overline\by$ and $\underline\by$ such that $-C_2\le
\underline\by\le\overline\by\le c_3$ for some $C_2>0$.

The other two cases can be proved similarly.
\end{pf}

\section{Small duration case revisited}\label{sec6}\label{sect-smallT}

In this and the next section, we shall argue that
the well-posedness of the dominating ODEs will lead to the desired
regular decoupling field. Our starting point will be the ``local
existence'' result for FBSDE, or more precisely, the well-posedness of
FBSDE (\ref{FBSDE}) over small time interval. We note that this
seemingly well-understood problem still contains many interesting
issues that have not been completely observed, especially in the case
when $\sigma $ depends on $z$ (i.e., $\sigma _3\neq0$), which we now
describe.

Let us first fix some constants $c_1, c_2>0$ such that
%
\begin{eqnarray}
\label{c2} c_1 c_2 < 1.
\end{eqnarray}
Set
$\tilde c_2 \stackrel{\triangle}{=}{c_2 + c_1^{-1}\over2}$, so that
$c_2 < \tilde c_2 < c_1^{-1}$.
Furthermore, recall $b_i$, $\sigma _i$, $f_i$, $i=1,2,3$ in~(\ref{tildephi}).
In what follows, all assumptions involving coefficients in (\ref{F1})
will be in the sense that they hold uniformly for all $\theta _j$, $j=1,2$.

Recall again that it is essential to have $\sigma _3 h\neq1$. We
shall establish the results for the cases (2-a)--(2-c) listed in
Section~\ref{sec5.2}. Our first result corresponds to case~(2-a) and Theorem
\ref{teo-nonlinear2}. We remark that the case $\sigma =\sigma
(t,x,y)$ satisfies case~(2-a) with arbitrary small $c_1>0$.

\begin{teo}
\label{teo-smallT1} Suppose that Assumption \ref{assum-Lipschitz} and
(\ref{c2}) are in force, and
assume that $\llvert  \sigma _3\rrvert  \le c_1$ and $\llvert  h\rrvert  \le c_2$.
Then there exists a constant $\delta >0$, which depends only on $c_1$,
$c_2$, and the
Lipschitz constants in Assumption \ref{assum-Lipschitz}, such that
whenever $T\le\delta $, it holds that:
\begin{longlist}[(iii)]
\item[(i)] the FBSDE (\ref{FBSDE}) has a unique solution $\Theta \in
\mathbb{L}^2$;

\item[(ii)] the ODEs in (\ref{ODE})
have solutions $\overline\by, \underline\by$ such that
%
\begin{eqnarray}
\label{ODEsmallT1} -\tilde c_2 \le\underline\by_t \le
\overline\by_t \le \tilde c_2\qquad\forall t\in[0,T];
\end{eqnarray}

\item[(iii)] there exists a random field $u$ such that, for all $t\in[0,T]$,
$Y_t = u(t, X_t)$ and
%
\begin{eqnarray}
\label{Yby1} \underline\by_t \le{u(t, x_1)-u(t, x_2) \over x_1-x_2} \le
\overline \by_t\qquad\mbox{for any } x_1 \neq
x_2.
\end{eqnarray}
\end{longlist}
\end{teo}

\begin{pf}
(i) follows directly from \cite{mybk} Theorem I.5.1. To
see (ii), we notice that $\underline F$ and $\overline F$ are
uniformly Lipschitz continuous in $y$ for $y\in[-\tilde c_2, \tilde
c_2]$ and denote by $L$ the uniform Lipschitz constant. We assume that
(i) holds for some $\delta >0$. Modifying $\delta $ if necessary we
may assume that
\begin{eqnarray*}
\biggl[\int_0^\delta e^{Lt}\,dt \biggr]
\Bigl[\sup_{\llvert  y\rrvert  \le\tilde c_2} \sup_{t\in[0,T]} \bigl[\bigl
\llvert \overline F(t, y)\bigr\rrvert +\bigl\llvert \underline F(t, y)\bigr\rrvert
\bigr] \Bigr] \le\tilde c_2-c_2.
\end{eqnarray*}
Now for any $T<\delta $, note that $\tilde\by^1 \stackrel{\triangle
}{=}- \tilde c_2$ and $\tilde\by^2 \stackrel{\triangle}{=}\tilde
c_2$ satisfy the following ODEs:
\begin{eqnarray*}
\tilde\by^1_t &=& -c_2 - [\tilde
c_2 - c_2] + \int_t^T
\bigl[\underline F\bigl(s, \tilde\by^1_s\bigr) -
\underline F(s, -\tilde c_2)\bigr]\,ds,
\\
\tilde\by^2_t &=& c_2 + [\tilde
c_2 - c_2] + \int_t^T
\bigl[\overline F\bigl(s, \tilde\by^2_s\bigr) - \overline
F(s, \tilde c_2)\bigr]\,ds.
\end{eqnarray*}
Following the arguments in Theorem \ref{teo-nonlinear2}, we prove (ii).

It remains to prove (iii). Let $\delta >0$ be small enough so that
both (i) and (ii) hold. For any $(t, x)$, denote the (unique) solution
to FBSDE (\ref{FBSDE}) starting from $(t,x)$ by $\Theta ^{t,x}$, and define
a random field $
u(t,x) \stackrel{\triangle}{=}Y^{t,x}_t$. The uniqueness of the
solution to FBSDE then leads to that $Y^{t,x}_s = u(s, X^{t,x}_s)$, for
all $s\in[t,T]$, $\mathbb{P}$-a.s. In particular, denoting $\Theta
_t=\Theta ^{0,x}_t$, we have $Y_t=u(t,X_t)$, $t\in[0,T]$.

Now let $x_1\neq x_2$ be given, and recall (\ref{tdX}) and (\ref
{linearFBSDE}).
Following standard arguments, see, for example, \cite{mybk} Theorem
I.5.1, for a smaller $\delta $ if necessary, one can easily see that
$\llvert  \nabla Y_t\rrvert  \le
\tilde c_2 \llvert  \nabla X_t\rrvert  $.
This also implies that
\begin{eqnarray*}
\llvert \nabla X_t\rrvert \le{1\over1-c_1\tilde c_2}\llvert \nabla
X_t - \sigma _3 \nabla Y_t\rrvert.
\end{eqnarray*}
Applying Theorem \ref{teo-linear-wellposed} we see that $\nabla X\neq
0$ and $\hat Y \stackrel{\triangle}{=}\nabla Y\slash\nabla X$
satisfies the BSDE
(\ref{hatYBSDE}) and (\ref{hatYproperty}). Then (\ref{Yby1}) follows
from Lemma \ref{lem-hatYbound}.
\end{pf}

Our next result corresponds to case (2-b) and Theorem \ref{teo-nonlinear3}.

\begin{teo}
\label{teo-smallT2} Suppose that Assumption \ref{assum-Lipschitz} and
(\ref{c2}) are in force, and assume that $\sigma _3$ and $h$
satisfy one of the conditions in (\ref{case21})--(\ref{case24}).
Then there exists a constant $\delta >0$, depending only on $c_1$,
$c_2$, and the
Lipschitz constants in Assumption \ref{assum-Lipschitz}, such that
when $T\le\delta $, all the results in Theorem \ref{teo-smallT1}
hold true, except that (\ref{ODEsmallT1}) should be replaced by the following:
\begin{eqnarray*}
\overline\by&\ge&\underline\by\ge\tilde c_2^{-1}\qquad\mbox{in cases (\ref{case21}) and (\ref{case22})}\quad\mbox{and}
\\
\underline \by&\le&\overline\by\le-\tilde c_2^{-1}\qquad\mbox{in cases (\ref{case23}) and (\ref{case24})}.
\end{eqnarray*}
%
\end{teo}

\begin{pf}
We shall argue that the assertions (i)--(iii) in Theorem
\ref{teo-smallT1} all remain true under the current assumptions.
Without loss of generality, we prove the result only for the case (\ref
{case21}). The other cases can be proved similarly.

We first assume (i) holds. Note that $c_2^{-1} \le h\le L$, where $L$
is the uniform Lipschitz constant in Assumption
\ref{assum-Lipschitz}. By similar arguments as those in Theorem
\ref{teo-smallT1}(ii), for $\delta $ small enough one can easily show
that the ODEs in (\ref{ODE})
have solutions $\overline\by, \underline\by$ such that
%
\begin{eqnarray}
\label{ODEsmallT2} \tilde c_2^{-1} \le\underline
\by_t \le\overline\by_t \le2L\qquad\mbox{for all } t\in[0,T].
\end{eqnarray}
This proves (ii). (iii) follows from (i) and similar arguments as those
in Theorem~\ref{teo-smallT1}(iii).

So it remains to prove (i). Our main idea is to reverse the roles of
forward and backward components and then apply Theorem \ref{teo-smallT1}.
To this end, we consider a simple transformation: $\tilde X \stackrel
{\triangle}{=}Y$ and $\tilde Y \stackrel{\triangle}{=}X$. In other
words, we define the coordinate change:
\begin{eqnarray*}
\lleft[ \matrix{ \tilde x
\cr
\tilde y}\rright] \stackrel{\triangle} {=} \lleft[ \matrix{
0 &1
\cr
1&0}\rright] \lleft[ \matrix{ x
\cr
y} \rright]\quad\mbox{and, correspondingly,}\quad
\tilde z \stackrel {\triangle} {=}\sigma (t,x,y,z).
\end{eqnarray*}
Note that, under (\ref{case21}), both functions $z\mapsto\sigma
(t,x,y,z)$ and $x\mapsto g(x)$ are invertible, that is, there exist
functions $\hat\sigma $ and $\hat g$ such that
%
\begin{eqnarray}
\label{inverse} \hat\sigma \bigl(t, x,y, \sigma (t,x,y, z)\bigr) = z,\qquad \hat
g \bigl(g(x)\bigr) = x.
\end{eqnarray}
Define
\begin{eqnarray*}
\tilde\sigma (t, \tilde\theta ) &\stackrel{\triangle} {=}&\hat \sigma (t, \tilde
y, \tilde x, \tilde z),\qquad \tilde g(\tilde x) \stackrel{\triangle} {=}\hat g(
\tilde x);
\\
\tilde b (t, \tilde\theta ) &\stackrel{\triangle} {=}&- f \bigl(t, \tilde y,
\tilde x, \tilde\sigma (t,\tilde\theta )\bigr),\qquad \tilde f (t, \tilde\theta )
\stackrel{\triangle} {=}- b\bigl(t, \tilde y, \tilde x, \tilde\sigma (t,\tilde
\theta )\bigr)
\end{eqnarray*}
and consider a new FBSDE:
%
\begin{eqnarray}
\label{tildeFBSDE} \cases{ \displaystyle \tilde X_t = \tilde x + \int
_0^t \tilde b(s, \tilde \Theta _s)
\,ds + \int_0^t \tilde\sigma (s,\tilde\Theta
_s) \,dB_s;
\cr
\displaystyle \tilde Y_t =
\tilde g(\tilde X_T) + \int_t^T
\tilde f(s, \tilde \Theta _s) \,ds - \int_t^T
\tilde Z_s \,dB_s,} \qquad t\in[0,T].
\end{eqnarray}

We now show that FBSDE (\ref{tildeFBSDE}) satisfies the conditions in
Theorem \ref{teo-smallT1}. First, by definition of inverse functions
and by (\ref{case21}), we have
\begin{eqnarray*}
\hat\sigma _1 + \hat\sigma _3 \sigma _1=0,
\qquad \hat\sigma _2 + \hat\sigma _3 \sigma
_2=0,\qquad \hat\sigma _3 \sigma _3=1\quad\mbox{and}\quad \hat h h =1,
\end{eqnarray*}
where $\hat\sigma _i, \hat h$ and more notation below are defined
in the spirit of (\ref{tildephi}) for the functions $\hat\sigma,
\hat g$. Note that
$\tilde\sigma _3 = \hat\sigma _3 = (\sigma _3)^{-1}$ and $\tilde h
= \hat h = h^{-1}$.
This implies that, by~(\ref{case21}),
%
\begin{eqnarray}
\label{sihbound1} L^{-1}\le\tilde\sigma _3\le
c_1,\qquad L^{-1}\le\tilde h \le c_2.
\end{eqnarray}
Next, since
\begin{eqnarray*}
\tilde b_1 = - f_2 - f_3 \tilde\sigma
_1 = - f_2 - f_3 \hat\sigma _2 =
- f_2 - f_3\sigma _2 (\sigma
_3)^{-1},
\end{eqnarray*}
we see that $\llvert  \tilde b_1\rrvert  \le C$. Similarly, $\llvert  \tilde\varphi _j\rrvert  \le C$
for $\varphi = b, \sigma, f$ and $j=1,2,3$.
Moreover, note that
\begin{eqnarray*}
\bigl\llvert \tilde g(0)\bigr\rrvert &=& \bigl\llvert \tilde g(0) - \tilde
g\bigl(g(0)\bigr)\bigr\rrvert \le L\bigl\llvert g(0)\bigr\rrvert;
\\
\bigl\llvert \tilde\sigma (t,0,0,0)\bigr\rrvert & =& \bigl\llvert \hat\sigma
(t,0,0,0)\bigr\rrvert = \bigl\llvert \hat\sigma (t,0,0,0) - \hat\sigma
\bigl(t,0,0,\sigma (t,0,0,0)\bigr)\bigr\rrvert
\\
&\le &C\bigl\llvert \sigma (t,0,0,0)
\bigr\rrvert;
\\
\bigl\llvert \tilde b(t,0,0,0)\bigr\rrvert &\le& \bigl\llvert f(t,0,0,0)\bigr
\rrvert +C\bigl\llvert \sigma (t,0,0,0)\bigr\rrvert,
\\
\bigl\llvert \tilde f(t,0,0,0)\bigr\rrvert &\le& \bigl\llvert b(t,0,0,0)\bigr\rrvert +C\bigl\llvert
\sigma (t,0,0,0)\bigr\rrvert.
\end{eqnarray*}
Thus (\ref{I0}) holds for FBSDE (\ref{tildeFBSDE}).

We can now apply Theorem \ref{teo-smallT1} to conclude that for some
$\delta >0$,
the FBSDE
(\ref{tildeFBSDE}) admits a unique solution $\tilde\Theta \in
\mathbb{L}^2$ for all $T\le\delta $, and $\tilde Y_t = \tilde u(t,
\tilde X_t)$
for some decoupling random field $\tilde u$. Moreover, by (\ref
{sihbound1}) and modifying the arguments in Theorem \ref{teo-smallT1}
slightly, we see that $\tilde u$ satisfies
%
\begin{eqnarray*}
{1\over2L} \le {\tilde u (t, \tilde x_1)- \tilde u(t, \tilde
x_2)\over
\tilde x_1 - \tilde x_2} \le\tilde c_2.
\end{eqnarray*}
Then $\tilde u(t,
\tilde x)$ has an inverse function $u(t,x)$ in terms of $x$. Now for
any $x$, let $\tilde x \stackrel{\triangle}{=}
u(0,x)$ and let $\tilde\Theta $ be the unique solution to FBSDE (\ref
{tildeFBSDE}) with initial value
$\tilde X_0 = \tilde x$. Then it is straightforward to check that
\[
X_t \stackrel{\triangle} {=}\tilde Y_t, \qquad
Y_t \stackrel {\triangle} {=}\tilde X_t, \qquad
Z_t \stackrel{\triangle} {=}\tilde \sigma (t, \tilde X_t,
\tilde Y_t, \tilde Z_t) %
\]
satisfy FBSDE (\ref{FBSDE}) with initial value $X_0=x$.

Finally, note that $\llvert  \tilde Z\rrvert   \le\llvert  \tilde\sigma (t,0,0,0)\rrvert   +
C[\llvert  \tilde X\rrvert  +\llvert  \tilde Y\rrvert  +\llvert  \tilde Z\rrvert  ]$, it is clear that $(X, Y, Z)\in
\mathbb{L}^2$. The proof is now complete.
\end{pf}

Our final result corresponds to case (2-c) and Theorem \ref{teo-nonlinear4}.

\begin{teo}
\label{teo-smallT3} Suppose that Assumption \ref{assum-Lipschitz} and
(\ref{c2})
are in force, and assume that $\sigma _3$ and $h$
satisfy one of the conditions in (\ref{case31})--(\ref{case34}).
Then there exists a constant $\delta >0$, depending only on $c_1$,
$c_2$, and the
Lipschitz constants in Assumption \ref{assum-Lipschitz}, such that
when $T\le\delta $, all the results in Theorem \ref{teo-smallT1}
hold true, except that (\ref{ODEsmallT1}) should be replaced by the following:
\begin{eqnarray*}
0 &\le&\underline\by\le\overline\by\le\tilde c_2,\qquad\mbox{in case of (\ref{case31})};
\\
\underline \by&\le&\overline\by\le \tilde c_2,\qquad\mbox{in case of (\ref{case32})};
\\
0&\ge& \overline\by\ge\underline\by\ge-\tilde c_2,\qquad\mbox{in case
of (\ref{case33})};
\\
\overline\by&\ge& \underline\by\ge -\tilde
c_2,\qquad\mbox{in case of (\ref{case34})}.
\end{eqnarray*}
%
\end{teo}

\begin{pf}
Again we consider only the case (\ref{case31}), and
the other cases can be argued similarly. Following similar arguments
as in Theorem \ref{teo-smallT2}, we shall only prove (i).

Slightly different from the proof of Theorem \ref{teo-smallT2} we
consider a slightly more complicated transformation:
$(\tilde x,\tilde y,\tilde z)\stackrel{\triangle}{=}\Phi
[\varepsilon ](x,y,z)$, where
%
\begin{eqnarray}
\label{tildexyz} \lleft[\matrix{\tilde x
\cr
\tilde y} \rright]\stackrel{\triangle}
{=} \lleft[\matrix{2\varepsilon &1
\cr
\varepsilon &1} \rright] \lleft[\matrix{ x
\cr
y} \rright], \qquad \tilde z \stackrel{\triangle} {=} \varepsilon \sigma (t, x,
y, z) + z.
\end{eqnarray}
Note that
%
\begin{eqnarray}
\label{case31-smallT} -L \le\sigma _3 \le c_1, \qquad 0\le h
\le c_2.
\end{eqnarray}
By choosing $\varepsilon >0$ small enough, we see that the mappings
\begin{eqnarray*}
z\mapsto\tilde z=\varepsilon \sigma (t, x,y,z) + z\quad\mbox{and}\quad  x\mapsto2
\varepsilon x + g(x)
\end{eqnarray*}
are both strictly increasing and thus both are invertible.
Denote the corresponding inverse functions by $\hat\sigma $ and $\hat
g$, respectively. Namely,
%
\begin{eqnarray}
\label{sig} \hat\sigma \bigl(t,x,y, \varepsilon \sigma (t,x,y,z)+z\bigr)=z,
\qquad \hat g\bigl(2\varepsilon x+g(x)\bigr)=x.
\end{eqnarray}
Furthermore, from (\ref{tildexyz}) we can solve
$ (x,y)=(\frac{\tilde x-\tilde y}{\varepsilon }, 2\tilde y-\tilde x)$,
the inverse transformation of $\Phi[\varepsilon ]$ is thus
\begin{eqnarray*}
(x,y,z)=\Psi[\varepsilon ](\tilde x,\tilde y,\tilde z)\stackrel {\triangle} {=}
\biggl(\frac{\tilde x-\tilde y}{\varepsilon }, 2\tilde y-\tilde x, \hat\sigma \biggl(t,
\frac{\tilde x-\tilde y}{\varepsilon
}, 2\tilde y-\tilde x, \tilde z \biggr) \biggr).
\end{eqnarray*}

We now consider the FBSDE (\ref{tildeFBSDE}) with the following new
coefficients:
%
\begin{eqnarray}
\label{tbf} \tilde b(t,\tilde x,\tilde y,\tilde z)&=&2\varepsilon b\bigl(t,
\Psi[\varepsilon ](\tilde x,\tilde y,\tilde z)\bigr)-f\bigl(t,\Psi[\varepsilon ](
\tilde x,\tilde y,\tilde z)\bigr),
\nonumber
\\
\tilde f(t,\tilde x,\tilde y,\tilde z)&=&-\varepsilon b\bigl(t,\Psi [\varepsilon
](\tilde x,\tilde y,\tilde z)\bigr)+f\bigl(t,\Psi[\varepsilon ](\tilde x,\tilde
y,\tilde z)\bigr),
\\
\tilde\sigma (t,\tilde x,\tilde y,\tilde z)&=&2\varepsilon \sigma \bigl(t,\Psi[
\varepsilon ](\tilde x,\tilde y,\tilde z)\bigr)+\hat\sigma \biggl(t,
\frac{\tilde x-\tilde y}{\varepsilon }, 2\tilde y-\tilde x, \tilde z\biggr),
\nonumber
\\
\tilde g(\tilde x) &=& \varepsilon \hat g(\tilde x) + g\bigl(\hat g(\tilde x)
\bigr).
\nonumber
\end{eqnarray}
Our idea is again to apply Theorem \ref{teo-smallT1}. Note that
$\hat\sigma _3 [\varepsilon \sigma _3 + 1] = 1$ and $\hat h
[2\varepsilon + h] = 1$, we have
\begin{eqnarray*}
\tilde\sigma _3 = 2\varepsilon \sigma _3 \hat\sigma
_3 + \hat \sigma _3 = {2\varepsilon \sigma _3 + 1\over\varepsilon \sigma _3
+ 1}, \qquad
\tilde h = \varepsilon \hat h + h \hat h = {\varepsilon
+ h\over2\varepsilon + h}.
\end{eqnarray*}
By (\ref{case31-smallT}) and for $\varepsilon >0$ small enough, we have
%
\begin{eqnarray}\label{tildesigmah}
0 &<& {1-2L \varepsilon \over1-\varepsilon L}\le\tilde\sigma _3 \le
{1+2c_1\varepsilon \over1+c_1\varepsilon }\stackrel{\triangle } {=}\overline c_1;
\nonumber\\[-8pt]\\[-8pt]\nonumber
0&<&{1\over2}\le\tilde h \le{\varepsilon +c_2\over2\varepsilon
+c_2}\stackrel{
\triangle} {=}\overline c_2.
\end{eqnarray}
Since $c_1c_2<1$, we obtain
%
\begin{eqnarray}
\label{barc1c2} \overline c_1 \overline c_2=
{1+2c_1\varepsilon \over
1+c_1\varepsilon } \cdot {\varepsilon +c_2\over2\varepsilon +c_2}<1.
\end{eqnarray}
Moreover, note that
$\hat\sigma _1 + \varepsilon \hat\sigma _3 \sigma _1 =0$ and
$\hat\sigma _2 + \varepsilon \hat\sigma _3\sigma _2=0$, we see that
$\hat\sigma _1 =
{-\varepsilon \sigma _1\over1+\varepsilon \sigma _3}$, $\hat
\sigma _2 = {-\varepsilon \sigma _2\over1+\varepsilon \sigma _3}$
are bounded and, therefore,
\begin{eqnarray*}
\tilde b_1 = 2\varepsilon \bigl[b_1 \varepsilon
^{-1} - b_2 + b_3\bigl[\hat\sigma
_1 \varepsilon ^{-1} - \hat\sigma _2\bigr]
\bigr] - \bigl[f_1 \varepsilon ^{-1} - f_2 +
f_3\bigl[\hat\sigma _1 \varepsilon ^{-1} -
\hat\sigma _2\bigr] \bigr]
\end{eqnarray*}
is
bounded. Similarly, one can check that all other coefficients are all
uniformly Lipschitz continuous and (\ref{I0}) also holds for FBSDE
(\ref{FBSDE}). Then we can apply Theorem \ref{teo-smallT1}, with
$c_1$, $c_2$ being replaced by $\overline c_1$, $\overline c_2$ here, to
conclude that (\ref{tildeFBSDE}) with coefficients given by
(\ref{tbf}) admits a unique solution $\tilde\Theta \in\mathbb
{L}^2$, for
$T\le\delta $ and $\delta $ small enough. Furthermore, by (\ref
{tildesigmah})
and following similar arguments as in Theorem~\ref{teo-smallT1}, it
holds that $\tilde Y_t = \tilde u(t, \tilde X_t)$ for some
decoupling random field $\tilde u$, which satisfies, for $\tilde
x_1 \neq\tilde x_2$, and $\overline c_3 \stackrel{\triangle
}{=}{\overline c_1^{-1}+\overline c_2\over
2}$,
\begin{eqnarray*}
{1\over4} \le{\tilde u(t, \tilde x_1) - \tilde u(t,
\tilde x_2)\over\tilde x_1-\tilde x_2}\le{\overline
c_3}.
\end{eqnarray*}
This then implies that $\tilde x\mapsto\tilde u(t, \tilde x)$ has an
inverse, denoted by $u(t,x)$.

Now for any $x$, let $\tilde x \stackrel{\triangle}{=}2\varepsilon
x+u(0,x)$ and $\tilde\Theta $ be the unique solution to FBSDE
(\ref{tildeFBSDE}) starting from $\tilde X_0 = \tilde x$. Then one can
easily check that
$\Theta:= \Psi[\varepsilon ](\tilde\Theta )$ satisfies all the
requirement.
\end{pf}

\section{Synthetic analysis}\label{sec7}
\label{sect-summary}

In this section, we summarize all the results proved in the previous
sections and
give a synthetic analysis for the solvability of FBSDE (\ref{FBSDE})
over an
arbitrary duration $[0,T]$, which in a sense could serve
as a \textit{User's Guide} for solving general FBSDEs. We should note that
all the cases listed
below cannot be
covered by the existing methods, therefore, they are all new.

\subsection{Linear case}\label{sec7.1}

We first consider the linear FBSDE (\ref{linearFBSDE}).
Bearing Remarks~\ref{sigma3h} and \ref{linearcase} in mind, then combining
Theorems \ref{teo-smallT1}, \ref{teo-smallT2} and \ref{teo-smallT3},
we have the following ``local'' well-posedness
result.
We note that since $\sigma $ is allowed to depend on $z$, and the
condition is both necessary and sufficient, this
result is already new.

\begin{teo}
\label{LFBSDE-smallduration}
Assume that the linear FBSDE (\ref{linearFBSDE}) has constant coefficients.
Then there exists a constant $\delta >0$, such that it is well-posed
on $[0,T]$, whenever
$T\le\delta $, if and only if
%
\begin{eqnarray}
\label{sigma3hles1} \sigma _3 h \neq1.
\end{eqnarray}
\end{teo}

\begin{rem}
\label{linearT} If the duration $T$ is arbitrarily given, then even
in the case when FBSDE is linear with constant coefficients the
necessary and sufficient conditions become slightly more
complicated. The reader should use Theorem \ref{teo-constant-F2} or
\ref{teo-constant-F3} as a benchmark.
\end{rem}

If the coefficients of FBSDE (\ref{linearFBSDE}) are random, then the
analysis becomes more involved. In fact, the
degree of difficulty is no less than that of general Lipschitz
coefficient case. We therefore do not discuss them separately.

\subsection{The case \texorpdfstring{$\sigma=\sigma(t,x,y)$}{$sigma=sigma(t,x,y)$}}\label{sec7.2}
We remark that the work \cite{Zhang} is a special case of the
following result.

\begin{teo}
\label{teo-wellposed1}
Assume all the conditions in Theorem \ref{teo-nonlinear1} hold, and
let $\underline \by, \overline \by$ be the bounded solutions of
ODEs (\ref{ODE}). Then:
\begin{longlist}[(iii)]
\item[(i)] FBSDE (\ref{FBSDE}) possesses a decoupling field $u$ satisfying
(\ref{Yby1}).

\item[(ii)] FBSDE (\ref{FBSDE}) admits a unique solution $\Theta \in\mathbb
{L}^2$, such
that
%
\begin{eqnarray}
\label{Norm} \llVert \Theta \rrVert _{\mathbb{L}^2}^2 \le C
\bigl[\llvert x\rrvert ^2+I_0^2\bigr].
\end{eqnarray}
Here, the constant $C>0$ depends only on $T$, the Lipschitz constant in
Assumption~\ref{assum-Lipschitz}, and the bound of $\underline \by,
\overline \by$.
\end{longlist}
\end{teo}

\begin{pf}
(i) First, applying Theorem \ref{teo-nonlinear1}, there
exists a constant $c_2>0$ such that
%
\begin{eqnarray}
\label{c2bound} -c_2 \le\underline \by_t \le\overline
\by_t \le c_2\qquad\mbox{for all } 0\le t\le T.
\end{eqnarray}
Notice that in this case $\sigma _3 =0$, thus we may set arbitrarily
small $c_1>0$ in Theorem~\ref{teo-smallT1}.

Let $\delta >0$ be the constant determined by $(c_1, c_2)$ in Theorem
\ref{teo-smallT1}, and $0=t_0<\cdots <t_n=T$ be a partition of
$[0,T]$ such that $t_i-t_{i-1}\le\delta $, $i=1,\ldots,n$. We first
consider FBSDE (\ref{FBSDE}) on $[t_{n-1}, t_n]$. Since $ \underline
\by_T \le h \le\overline \by_T$, we see\vspace*{1pt} that the Lipschitz constant
of the terminal condition $g$ is less than $ c_2$, then by Theorem \ref
{teo-smallT1} there exists a random field $u(t,x)$ for $t\in[t_{n-1},
t_n]$ such that (\ref{Yby1}) holds for all $t\in[t_{n-1}, t_n]$. In
particular, the estimate (\ref{Yby1}) at $t_{n-1}$ and (\ref
{c2bound}) imply that $c_2$ is also a Lipschitz constant of $u(t_{n-1},
\cdot )$. Next, consider FBSDE (\ref{FBSDE}) on $[t_{n-2}, t_{n-1}]$
with terminal condition $u(t_{n-1}, \cdot )$. Applying Theorem \ref
{teo-smallT1} again, we find $u$ on $[t_{n-2}, t_{n-1}]$ such that
(\ref{Yby1}) holds for $t\in[t_{n-2}, t_{n-1}]$. Repeating this
procedure backwardly finitely many times, we extend the random field
$u$ to the whole interval $[0, T]$. Clearly, it is a decoupling field
satisfying (\ref{Yby1}).

(ii) We first note that the above $n$ is fixed. Since $u$ is uniformly
Lipschitz continuous in $x$, applying Theorem \ref{teo-smallT1} on
each interval $[t_i, t_{i+1}]$ with initial value $X_{t_i}=0$, we see
that there exists a constant $C$ such that
\begin{eqnarray*}
\mathbb{E}\bigl\{\bigl\llvert u(t_i,0)\bigr\rrvert ^2
\bigr\} = \mathbb{E}\bigl\{\bigl\llvert Y^{t_i,0}_{t_i}\bigr
\rrvert ^2\bigr\} \le C\mathbb{E}\bigl\{ \bigl\llvert
u(t_{i+1},0)\bigr\rrvert ^2\bigr\} + CI_0^2.
\end{eqnarray*}
Note that $u(t_n,0) = g(0)$, we see that, for a larger $C$, $\max_{0\le i\le n} \mathbb{E}\{\llvert  u(t_i,0)\rrvert  ^2\} \le CI_0^2$.

Next, by Theorem \ref{teo-decoupling} FBSDE (\ref{FBSDE}) admits a
unique global solution $\Theta $. Applying Theorem \ref{teo-smallT1}
on each interval $[t_i, t_{i+1}]$ again, we obtain
%
\begin{eqnarray}\label{Normi}
&& \mathbb{E} \biggl\{\sup_{t_i\le t\le t_{i+1}}\bigl[\llvert
X_t\rrvert ^2 + \llvert Y_t\rrvert
^2\bigr] + \int_{t_i}^{t_{i+1}}Z_t^2\,dt
\biggr\}
\nonumber\\[-8pt]\\[-8pt]\nonumber
&&\qquad  \le C\mathbb{E} \bigl\{\llvert X_{t_i}\rrvert ^2 +
\bigl\llvert u(t_{i+1},0)\bigr\rrvert ^2 \bigr\} +
CI_0^2.
\end{eqnarray}
This implies that
\begin{eqnarray*}
\mathbb{E}\bigl\{\llvert X_{t_{i+1}}\rrvert ^2\bigr\} \le C
\mathbb{E} \bigl\{\llvert X_{t_i}\rrvert ^2 + \bigl\llvert
u(t_{i+1},0)\bigr\rrvert ^2 \bigr\} + CI_0^2
\le C\mathbb{E}\bigl\{\llvert X_{t_i}\rrvert ^2\bigr\} +
CI_0^2,
\end{eqnarray*}
thus $ \max_{i} \mathbb{E}\{\llvert  X_{t_i}\rrvert  ^2\} \le C[\llvert  x\rrvert  ^2+I_0^2]$.
Plugging into (\ref{Normi}) and summing over $i$, we obtain (\ref{Norm}).
\end{pf}

In Table \ref{tab1} below we list a few classes
of FBSDEs whose coefficients $(b, \sigma,
f)$ satisfy condition (\ref{nonlinear-small0}), and thus are
well-posed for arbitrary $T$ under standard Lipschitz conditions. We
note that all coefficients are allowed to be
random, and
$Y_T=g(X_T)$.
%

\begin{table}[t]%
\tabcolsep=0pt
\tablewidth=280pt
\caption{Cases satisfying (\protect\ref{nonlinear-small0})}\label{tab1}
\begin{tabular*}{\tablewidth}{@{\extracolsep{\fill}}@{}lccc@{}}
\hline
\textbf{Assumption} & $\bolds{b}$& $\bolds{\sigma}$& $\bolds{f}$\\
\hline
& $b(t,x,z)$   &   $\sigma(t)$   &   $f(t,x,y,z)$   \\
&   $b(t,x)$ &  $\sigma(t,x,y)$ &     $f(t,x,y)$\\
$\sigma _2 b_3 \le0$, $\beta_t\ge c$ &   $b(t,x, z)$ & $\sigma(t,\beta_t x + y)$  &     $f(t,x, y)$\\
$\sigma _2 b_3 \le0$, $\beta_t\ge c$ &   $b(t,x, z)$ &  $\sigma(t, y)$ &$f_0(t,x, y) + \beta_t b(t,x, z)$\\
\hline
\end{tabular*}
\end{table}

\subsection{The general case \texorpdfstring{$\sigma=\sigma(t,x,y,z)$}{$sigma=sigma(t,x,y,z)$}}\label{sec7.3}
We now turn to the general case. We assume that the standing
Assumption \ref{assum-Lipschitz}, (\ref{c}), and one of the
assumptions (\ref{case1}), (\ref{case21})--(\ref{case24}) and
(\ref{case31})--(\ref{case34}) hold. For the convenience of the
reader, we tabulate these conditions so that
the nature of these assumptions are more explicit. Let $\varepsilon
>0$ be given as that
in Theorems \ref{teo-nonlinear2}, \ref{teo-nonlinear3}, \ref
{teo-nonlinear4} and $\alpha _3\stackrel{\triangle}{=}b_2-\frac
{b_3\sigma _2}{\sigma _3}$.

\begin{longlist}[\textit{Case} III.]
\item[\textit{Case} I.] $\llvert  \sigma _3\rrvert  \le c_1$, $\llvert  h\rrvert  \le c_2$; and
$ \overline F(t, c_3) \le\varepsilon $, $\underline F(t, -c_3) \ge
-\varepsilon $.

\item[\textit{Case} II.] $\llvert  \sigma _3\rrvert  \ge c_1^{-1}$, $\llvert  h\rrvert  \ge
c_2^{-1}$, and both of them keep the same sign (see Table~\ref{tab2}).

\begin{table}[b]
\tabcolsep=0pt
\tablewidth=280pt
\caption{$\sigma_3\neq 0$, Case~\textup{II}}\label{tab2}
\begin{tabular*}{\tablewidth}{@{\extracolsep{\fill}}@{}lcc@{}}
\hline
& $\bolds{h \ge c_2^{-1}}$ & $\bolds{h \le-c_2^{-1}}$\\
\hline
$\sigma _3 \ge c_1^{-1}$& $\underline F(t, c_3^{-1}) \ge-\varepsilon$, $\alpha _3 \le\varepsilon $ & $\overline F(t, c_3^{-1}) \le\varepsilon $, $\alpha _3 \ge -\varepsilon $\\
$\sigma _3 \le-c_1^{-1}$& $\underline F(t, c_3^{-1}) \ge -\varepsilon $, $\alpha _3 \le\varepsilon $ & $\overline F(t, c_3^{-1}) \le\varepsilon $, $\alpha _3 \ge -\varepsilon $\\
\hline
\end{tabular*}
\end{table}

\item[\textit{Case} III.] $\sigma _3 h \le c_1c_2$, and one of them
keeps the same sign (see Table~\ref{tab3}).

Our main result is the following.

\begin{teo}
\label{teo-wellposed}
Suppose that Assumption
\ref{assum-Lipschitz} and (\ref{c}) are in force, and for $c_1$,
$c_2$, $c_3$ in
(\ref{c}), either one of the conditions listed in cases \textup{I--III} holds.
Then:
\begin{longlist}[(iii)]
\item[(i)] FBSDE (\ref{FBSDE}) possesses a decoupling field $u$ such that
${u(t,x_1)-u(t,x_2)\over x_1-x_2}$ satisfies the corresponding property
of $h$ with $c_2$ being replaced by $c_3$.

\item[(ii)] FBSDE (\ref{FBSDE}) admits a unique solution $\Theta \in\mathbb
{L}^2$, and there exists a constant $C>0$, depending only on $T$, the
Lipschitz constant in Assumption \ref{assum-Lipschitz}, and $c_1, c_2,
c_3$, such that (\ref{Norm}) holds.
\end{longlist}
\end{teo}
\end{longlist}

\begin{pf}
The proof is similar to that of Theorem \ref{teo-wellposed1} and is thus omitted. However, we emphasize that when
one applies Theorems \ref{teo-smallT1}, \ref{teo-smallT2} or \ref
{teo-smallT3}, the constant $\delta $ should be determined by $c_1,
c_3$, not by $c_1, c_2$.
\end{pf}

The following special case deserves special attention.

\begin{cor}
\label{cor-wellposed} Assume that Assumption \ref{assum-Lipschitz}
hold. If the coefficients
in the variational FBSDE (\ref{variFBSDE}), defined by
(\ref{tildephi}), satisfy either
%
\begin{eqnarray}
\label{case41} \sigma _3 \ge0, \qquad h\le0, \qquad f_1
\le0,\qquad b_2 - {b_3
\sigma _2\over\sigma _3} \ge0
\end{eqnarray}
or
%
\begin{eqnarray}
\label{case42} \sigma _3 \le0, \qquad h\ge0, \qquad f_1
\ge0,\qquad b_2 - {b_3 \sigma _2\over
\sigma _3} \le0;
\end{eqnarray}
then the FBSDE (\ref{FBSDE}) is well-posed over arbitrary duration $[0,T]$.
\end{cor}

\begin{pf}
We assume (\ref{case41}) holds. Let $c_1$ be the Lipschitz constant
of $\sigma $ with respect to $z$, and let $0< c_2 <c_3<\delta $ for
some $\delta $ small enough. One can easily check that (\ref
{case32}) holds.
\end{pf}

\subsection{Comparison to the existing methods}\label{sec7.4}

We now compare our conditions to
those of the three well-known existing methods.

\begin{table}
\tabcolsep=0pt
\tablewidth=280pt
\caption{$\sigma_3\neq 0$, Case~\textup{III}}\label{tab3}
\begin{tabular*}{\tablewidth}{@{\extracolsep{\fill}}@{}lcccc@{}}
\hline
&  $\bolds{h \ge0}$ &   $\bolds{h\le0}$ &   $\bolds{\sigma _3 \ge0}$ &   $\bolds{\sigma _3\le0}$\\
\hline
$\sigma _3 \le c_1$, & $\overline F(t, c_3) \le\varepsilon $, & & $\overline F(t, c_3) \le\varepsilon $,&\\
$h\le c_2$ & $f_1\ge0$ && $\alpha _3 \ge-\varepsilon $&
\\[6pt]
$\sigma _3 \ge-c_1$, & & $\underline F(t, -c_3) \ge-\varepsilon $,& &$\underline F(t, -c_3) \ge-\varepsilon $, \\
$h\ge-c_2$ & &     $f_1\le0$ &&    $\alpha _3 \le\varepsilon $\\
\hline
\end{tabular*}
\end{table}

\begin{longlist}
\item[1. \textit{Method of Contraction Mapping.}]
It has been understood that the
fundamental assumptions for this method are
$\llvert  \sigma _3 g_1\rrvert  < 1$ and that $T$ is small enough (see, e.g., \cite{mybk}, Theorem I.5.1). In fact, \cite{mybk}, Example I.5.2,
shows that the FBSDE could be unsolvable if $\sigma _3g_1=1$.
Therefore, Theorem \ref{LFBSDE-smallduration} in this paper
indeed presents the sharpest result in the linear case.

For the general case, we note that in Antonelli \cite{fabio} $\sigma
_3=0$. To compare with the work of Pardoux and Tang \cite{PT}, we recall
(\ref{F}). Then it is easy to see that
in \cite{PT} it is essentially assumed, besides $\sigma _3$ and $h$
satisfying condition (\ref{case1}), that one of the following
conditions holds:
\begin{longlist}[(ii)]
\item[(i)] either $b_2, b_3, \sigma _2, \sigma _3$ or $f_1, h$ are small
(``\textit{weak coupling}'');

\item[(ii)] either $b_1$ or $f_2$ is very negative (``\textit{strong monotone}'').

But
for fixed $T$, (i) implies that the coefficients
of $y^2$ and $y^3$ is small enough, and thus the ODEs (\ref{ODE}) has
desired solutions on $[0, T]$, and (ii)
implies that the coefficient of $y$ is very negative, which ensures
that the solution to ODEs (\ref{ODE}) does
not blow up before $T$.
\end{longlist}

\item[2. \textit{Method of Continuation.}] The ``monotonicity condition'' in
Hu and Peng \cite{hupeng},
Peng and Wu \cite{PW}, Yong \cite{yong1} states
%
\begin{eqnarray}
\label{mono} \Delta b\Delta y + \Delta \sigma \Delta z - \Delta f \Delta x
&\ge& \beta\bigl[\llvert \Delta x\rrvert ^2 + \llvert \Delta y\rrvert
^2 + \llvert \Delta z\rrvert ^2\bigr],
\nonumber\\[-8pt]\\[-8pt]\nonumber
\Delta g \Delta x &\le& - \beta\llvert \Delta x\rrvert ^2,
\end{eqnarray}
for some constant $\beta>0$. By some simple analysis, one sees
immediately that (\ref{mono}) implies
\begin{eqnarray*}
b_2\ge\beta, \qquad \sigma _3 \ge\beta, \qquad
f_1 \le-\beta\le 0,\qquad h\le-\beta\le0.
\end{eqnarray*}
Moreover, by setting $\Delta x = 0$, we see that
\begin{eqnarray*}
b_2 \llvert \Delta y\rrvert ^2 +\sigma _3
\llvert \Delta z\rrvert ^2 + (b_3 + \sigma _2)
\Delta y\Delta z\ge0\qquad\mbox{for any } \Delta y,\Delta z.
\end{eqnarray*}
Then it must hold that $(b_3 + \sigma _2)^2 - 4 b_2 \sigma _3 \le0$,
and thus $b_2\sigma _3 \ge{1\over4} (b_3 + \sigma _2)^2 \ge
b_3\sigma _2$. These lead exactly to (\ref{case41}), and thus the
FBSDE is well-posed. Clearly, the monotonicity condition
can be easily further weakened in our framework.
\end{longlist}

\begin{longlist}
\item[3. \textit{Four Step Scheme.}] We should note that our solvability
conditions\break (\ref{case1}), (\ref{case21})--(\ref{case24}), (\ref
{case31})--(\ref{case34}) do not cover the
results in \cite{mpy} and \cite{Delarue}. This is because the
generality of the FBSDE that we are pursuing in this paper, especially
the non-Markovian structure (i.e., random coefficients) and the
possible degeneracy of $\sigma $, essentially
inhibits us from taking advantage of the special features of
nondegenerate PDEs. We nevertheless observe that in both \cite{mpy}
and \cite{Delarue}, the solution of the PDE, which serves as a
deterministic decoupling function, is indeed uniformly Lipschitz
continuous, and thus
falls into the framework of Theorem \ref{teo-decoupling}. In fact, our
definition of regular decoupling fields is strongly motivated by these works.
\end{longlist}

\subsection{Regarding examples (\texorpdfstring{\protect\ref{eg1}}{1.2}) and (\texorpdfstring{\protect\ref{eg2}}{1.3})}\label{sec7.5}
%
We now return to the two examples (\ref{eg1}) and (\ref{eg2})
mentioned in the \hyperref[sec1]{Introduction}. Note that in (\ref{eg1})
we actually have
$F(h) = 0$ and $b_2 - {b_3\sigma _2\over\sigma _3}=0$. Then,\vspace*{1.5pt} for
$\sigma \neq0, 1$, either (i) or (ii) of Theorem \ref{teo-constant-F3} will hold, and thus the FBSDE is well-posed. Since
the equation is trivial for $\sigma =0$, we can thus
conclude that
\textit{the FBSDE} (\ref{eg1}) \textit{is well-posed if and only if $\sigma \neq1$}.

We now turn attention to example (\ref{eg2}). To understand the
problem, we briefly describe its origin (see
\cite{CZ} for more details). Consider the following FBSDE:
%
\begin{eqnarray}
\label{FBSDE-CZ} dX_t &=& \sigma (t, X_t,
Y_t) \,dB_t,\qquad dY_t = f(t, X_t,
Y_t)\,dt -Z_t \,dB_t;
\nonumber\\[-8pt]\\[-8pt]\nonumber
X_0&=&x, \qquad Y_T=g(X_T),
\end{eqnarray}
where the coefficients are all deterministic. The purpose is to find a
Monte Carlo method for the numerical solution,
without using PDEs.
Following the idea of ``method of optimal control'' (cf. \cite{mybk}),
one can consider (\ref{FBSDE-CZ}) as a controlled
diffusion starting from $(x,y)$, and try to find the ``control'' $(y,
Z)\in\mathbb{R}\times L^2_{\mathbb{F}}([0,T])$ so that
%
\begin{eqnarray*}
0=\inf_{y, Z} V(x, y;Z) \stackrel{\triangle} {=}
\inf_{y,Z}{1\over
2}\mathbb{E} \bigl[\bigl\llvert
Y^{x,y,Z}_T - g\bigl(X^{x, y,Z}_T\bigr)\bigr
\rrvert ^2 \bigr].
\end{eqnarray*}
Since the existence of the optimal control is known (as the FBSDE is
solvable), the main task here is to numerically
compute the optimal control and trajectory. We proceed iteratively:
given some initial control $(y_0, Z^0)$ and we
find the approximating sequence $(y_n, Z^n)$ that converges to the true
solution $(Y_0, Z)$ of the FBSDE (\ref{FBSDE-CZ}). The so-called
``steepest descent method'' proposed in \cite{CZ} suggests that at
each step one should set
$ ( y_n, Z^n):= (y_{n-1}, Z^{n-1})-\lambda (\overline Y^n_0, \overline Z^n)$
for some small constant $\lambda >0$,
where $(\overline Y^n, \overline Z^n, \tilde Y^n,\break \tilde Z^n)$ solves a certain
BSDE which can be rewritten as
%
\begin{eqnarray}
\label{adjoint2} \overline Y^n_t &=&\overline
Y_0 + \int_0^t
\bigl[f_y \overline Y^n_s + \sigma
_y \tilde Z^n_s \bigr]\,ds + \int
_0^t \overline Z^n_s
\,dB_s;
\nonumber\\[-8pt]\\[-8pt]\nonumber
\tilde Y^n_t &=& g_x \overline
Y^n_T + \int_t^T
\bigl[f_x \overline Y^n_s + \sigma
_x \tilde Z^n_s \bigr]\,ds - \int
_t^T \tilde Z^n_s
\,dB_s.
\end{eqnarray}
If we
view $\overline Z^n$ as a given random coefficient, $\overline Y^n$ the
forward component, and $(\tilde Y^n, \tilde Z^n)$ the backward one,
then equations (\ref{adjoint2}) is an FBSDE same as (\ref{eg2}).
This FBSDE cannot be covered by any existing method, but it
satisfies condition (\ref{nonlinear-small0}), and thus falls into our
framework.
Furthermore applying Corollary
\ref{cor-Lp} below we can derive an important estimate in \cite{CZ}.
We refer the interested reader to \cite{CZ} for details.
%
%


\section{Properties of the solution}\label{sec8}
\label{sect-property}

In this section, we establish some further properties of the solution
to the
FBSDE (\ref{FBSDE}). These will include a stability result, an
$\mathbb{L}^p$-estimate for $p>2$, and a comparison theorem for FBSDE.

We first prove the stability result.

\begin{teo}[(Stability)] Assume both $(b, \sigma, f, g)$ and $(\tilde
b, \tilde\sigma,
\tilde f, \tilde g)$ satisfy the same conditions (i.e., they belong
to the same case) in Theorem~\ref{teo-wellposed} (or Theorem~\ref{teo-wellposed1}). Let $u, \tilde
u$ be the
corresponding random fields and, for any $(t,x)$, $\Theta ^{t,x}$ and
$\tilde\Theta ^{t,x}$ the solutions to the corresponding FBSDEs.
For $\varphi = b, \sigma, f, g$, denote $\Delta \varphi \stackrel
{\triangle}{=}\tilde\varphi -\varphi $. Then
%
\begin{eqnarray}
&&\bigl\llVert \tilde\Theta ^{0,\tilde x} - \Theta
^{0,x}\bigr\rrVert ^2_{\mathbb{L}^2}\nonumber
\\
\label{DThest}  &&\qquad \le C\mathbb{E} \biggl\{\llvert \tilde x-x\rrvert ^2 + \bigl
\llvert \Delta g\bigl(X^{0,x}_T\bigr)\bigr\rrvert
^2
\\
&&\hspace*{54pt}{} + \biggl(\int_0^T \bigl[\llvert
\Delta b\rrvert +\llvert \Delta f\rrvert \bigr]\bigl(t,\Theta
^{0,x}_t\bigr)\,dt \biggr)^2 + \int
_0^T \llvert \Delta \sigma \rrvert
^2\bigl(t,\Theta ^{0,x}_t\bigr)\,dt \biggr\},
\nonumber\hspace*{-30pt}
\\
&&\bigl\llvert \tilde u(t,x)-u(t,x)\bigr\rrvert ^2\nonumber
\\
\label{Duest} &&\qquad  \le C\mathbb{E}_t \biggl\{\bigl\llvert \Delta g
\bigl(X^{t,x}_T\bigr)\bigr\rrvert ^2+ \biggl(
\int_t^T \bigl[\llvert \Delta b\rrvert +\llvert
\Delta f\rrvert \bigr]\bigl(s,\Theta ^{t,x}_s\bigr)\,ds
\biggr)^2
\\
&&\hspace*{165pt}{}  + \int_t^T \llvert \Delta
\sigma \rrvert ^2\bigl(s,\Theta ^{t,x}_s\bigr)\,ds
\biggr\} \qquad \mbox {a.s.}
\nonumber
\end{eqnarray}
\end{teo}

\begin{pf}
Note that
$\tilde u(t,x)-u(t,x) = \tilde Y^{t,x}_t - Y^{t,x}_t$, and
Consider the FBSDEs on $[t,T]$ and replace $\mathbb{E}$ with $\mathbb
{E}_t$, (\ref{Duest}) follows directly from (\ref{DThest}).

To show (\ref{DThest}), denote $\Delta \Theta \stackrel{\triangle
}{=}\tilde\Theta ^{0,\tilde x} - \Theta ^{0,x}$ and $\Delta
x\stackrel{\triangle}{=}\tilde x-x$. Then
\begin{eqnarray*}
\Delta X_t &=& \Delta x + \int_0^t
\bigl[\tilde b_1 \Delta X_s + \tilde b_2
\Delta Y_s + \tilde b_3 \Delta Z_s + \Delta
b\bigl(s,\Theta ^{0,x}_s\bigr)\bigr]\,ds
\\
&&{} + \int_0^t \bigl[\tilde\sigma
_1 \Delta X_s + \tilde\sigma _2 \Delta
Y_s +\tilde\sigma _3 \Delta Z_s + \Delta
\sigma \bigl(s,\Theta ^{0,x}_s\bigr)\bigr]\,dB_s;
\\
\Delta Y_t &=& \tilde h \Delta X_T + \Delta g
\bigl(X^{0,x}_T\bigr)
\\
&&{} + \int_t^T
\bigl[\tilde f_1 \Delta X_s + \tilde f_2
\Delta Y_s + \tilde f_3 \Delta Z_s + \Delta
f\bigl(s,\Theta ^{0,x}_s\bigr)\bigr]\,ds - \int
_t^T \Delta Z_s \,dB_s.
\end{eqnarray*}
Here, the notation $\tilde b_1$, etc., are defined similar to (\ref
{tildephi}). One can easily check that the above linear FBSDE (with
solution $\Delta \Theta $) satisfies the corresponding conditions in
Theorem \ref{teo-wellposed} (or Theorem \ref{teo-wellposed1}). Then
applying the theorem we obtain the estimate immediately.
\end{pf}

We next establish the $L^p$-estimates for some $p>2$. First, following
Karatzas and Shreve \cite{KS1} (cases 2 and 4, page~164), one can easily
prove the following lemma.

\begin{lem}
\label{lem-Lp}
For any $p\ge2$ and $Z\in L^{2,p}$, that is, $E [ (\int_0^T
\llvert  Z_t\rrvert  ^2\,dt )^{p/2} ]<\infty$, we have
%
\begin{eqnarray}
\label{MLp} \bigl\llvert \psi_1(p)\bigr\rrvert ^{-p}
E \biggl[\biggl\llvert \int_0^t
Z_s \,dB_s\biggr\rrvert ^p \biggr] &\le& E
\biggl[ \biggl(\int_0^t \llvert
Z_s\rrvert ^2\,ds \biggr)^{p/2} \biggr]
\nonumber\\[-8pt]\\[-8pt]\nonumber
&\le&\bigl
\llvert \psi _2(p)\bigr\rrvert ^{p} E \biggl[\biggl\llvert
\int_0^t Z_s \,dB_s
\biggr\rrvert ^p \biggr],
\end{eqnarray}
where
%
\begin{eqnarray}
\label{psi} \psi_1(p) &\stackrel{\triangle} {=}& 2^{-{1/ p}}p^{1/2}
\biggl({2^{p/2}-2\over p-2} \biggr)^{{1/2}-{1/ p}},
\nonumber\\[-8pt]\\[-8pt]\nonumber
\psi _2(p) &\stackrel{\triangle} {=}& \biggl({p-1\over2} \biggr)^{1/ p}
p^{1/2} \biggl[2^{p/2} + {2^{{p/2}}-2/(p-2)}
\biggr]^{{1/2}-{1/p}}.
\end{eqnarray}
Moreover, for $i=1, 2$,
$\psi_i$ is continuous, strictly increasing on $[2, \infty)$ and 
$\psi_i(2) = 1$, $\psi_i(\infty) = \infty$.
\end{lem}

We now give the $L^p$-estimate of the solutions.

\begin{teo}[($L^p$-estimates)]
\label{teo-Lp} Let $(b,\sigma,f,g)$ satisfy the conditions in
Theorem \ref{teo-wellposed}. Assume
%
\begin{eqnarray}
\label{p} 2\le p < \psi^{-1}\biggl({1\over c_1c_3}
\biggr),
\end{eqnarray}
where $\psi\stackrel{\triangle}{=}\psi_1\psi_2$ and $\psi^{-1}$
denote the inverse function of $\psi$; and
%
\begin{eqnarray}
\label{Ip} I_p^p &\stackrel{\triangle} {=}&
\mathbb{E} \biggl\{ \biggl(\int_0^T \bigl[
\llvert b\rrvert +\llvert f\rrvert \bigr](t,0,0,0)\,dt \biggr)^p
\nonumber\\[-8pt]\\[-8pt]\nonumber
&&\hspace*{12pt}{} +
\biggl(\int_0^T \llvert \sigma \rrvert
^2(t,0,0,0)\,dt \biggr)^{p/2}+ \bigl\llvert g(0)\bigr\rrvert
^p \biggr\}<\infty.
\end{eqnarray}
Then the unique solution $\Theta $ of FBSDE (\ref{FBSDE}) is in $L^p$
and satisfies
%
\begin{eqnarray}
\label{Lp} \llVert \Theta \rrVert _{L^p}\le C_p \bigl[
\llvert x\rrvert +I_p \bigr].
\end{eqnarray}
Consequently, the corresponding
characteristic BSDE (\ref{hatYBSDE}) has a unique solution
$(\hat Y, \hat Z)$
satisfying (\ref{hatYproperty}) and (\ref{hatZproperty}).
\end{teo}

\begin{pf}
By Theorem \ref{teo-wellposed} and following its
arguments, we may assume $p>2$ and shall only prove the theorem under
(\ref{case1}) and for $T\le\delta $, where $\delta $ is a constant
which depends on $c_1, c_3$, the Lipschitz constants, and $p$ and will
be specified later. Moreover, by using the standard stopping arguments,
we can assume without loss of generality that
%
\begin{eqnarray}
\label{ThLp0} \llVert \Theta \rrVert _{w,p}^p\stackrel{
\triangle} {=}\mathbb{E} \biggl[\int_0^T \bigl[
\llvert X_t\rrvert ^{p}+\llvert Y_t\rrvert
^{p}\bigr]\,dt+\biggl(\int_0^T
\llvert Z_t\rrvert ^{2}\,dt\biggr)^{{p}/{2}} \biggr]<
\infty.
\end{eqnarray}

For any $0<\varepsilon \le1$ and $a, b>0$, note that $(a+b)^p \le
C_{p,\varepsilon } a^p + (1+\varepsilon ) b^p$, for some generic
constant $C_{p,\varepsilon }\ge1$ which may depend on $p$ and
$\varepsilon $. Then, for any $0\le t\le T\le\delta $,
we have [denoting $\varphi _s=\varphi (s,\Theta _s)$, $\varphi
=b,\sigma $, for simplicity]
\begin{eqnarray*}
\mathbb{E}\bigl[\llvert X_t\rrvert ^p\bigr] &\le&
C_{p,\varepsilon }\mathbb{E} \biggl[\llvert x\rrvert ^p+ \biggl(\int
_0^t \llvert b_s\rrvert \,ds
\biggr)^p \biggr] + (1+\varepsilon )\mathbb {E} \biggl[ \biggl\llvert
\int_0^t \sigma _s\,dB_s
\biggr\rrvert ^p \biggr]
\\
&\le& C_{p,\varepsilon }\mathbb{E} \biggl[\llvert x\rrvert ^p+
\biggl(\int_0^t \llvert b_s\rrvert
\,ds \biggr)^p \biggr]
\\
&&{}  + (1+\varepsilon )\psi_1(p)^p
\mathbb{E} \biggl[ \biggl(\int_0^t \llvert
\sigma _s\rrvert ^2\,ds \biggr)^{p/2} \biggr],
\end{eqnarray*}
where the second inequality thanks to Lemma \ref{lem-Lp}. Note that
\begin{eqnarray*}
&& \biggl[\int_0^t \llvert b_s
\rrvert \,ds \biggr]^p
\\
&&\qquad \le  C_p \biggl\{\int
_0^T \bigl[\bigl\llvert b(s, 0)\bigr\rrvert +
\llvert X_s\rrvert +\llvert Y_s\rrvert +\llvert
Z_s\rrvert \bigr]\,ds \biggr\}^p
\\
&&\qquad \le C_p \biggl\{ \biggl[\int_0^T
\bigl\llvert b(s, 0)\bigr\rrvert \,ds \biggr]^p
\\
&&\quad\qquad{} + T^{p-1}\int
_0^T \bigl[\llvert X_s\rrvert
^p+\llvert Y_s\rrvert ^p\bigr]\,ds +
T^{p/2} \biggl[\int_0^T\llvert
Z_s\rrvert ^2\,ds \biggr]^{p/2} \biggr\},
\\
&& \biggl [\int_0^t \llvert \sigma
_s\rrvert ^2\,ds \biggr]^{p/2}
\\
&&\qquad \le \biggl(\int
_0^T \bigl[C_{\varepsilon }\bigl[\bigl\llvert
\sigma (s, 0)\bigr\rrvert ^2+ \llvert X_s\rrvert
^2+\llvert Y_s\rrvert ^2\bigr] + (1+
\varepsilon )c_1^2\llvert Z_s\rrvert
^2 \bigr]\,ds \biggr)^{p/2}
\\
&&\qquad \le C_{p,\varepsilon } \biggl\{\int_0^T
C_{\varepsilon
}\bigl[\bigl\llvert \sigma (s, 0)\bigr\rrvert ^2+
\llvert X_s\rrvert ^2+\llvert Y_s\rrvert
^2\bigr]\,ds \biggr\}^{p/2}
\\
&&\quad\qquad{}  + (1+\varepsilon ) \biggl[\int
_0^T (1+\varepsilon )c_1^2
\llvert Z_s\rrvert ^2\,ds \biggr]^{p/2}
\\
&&\qquad \le C_{p,\varepsilon } \biggl\{\int_0^T
\bigl[\bigl\llvert \sigma (s, 0)\bigr\rrvert ^2+ \llvert
X_s\rrvert ^2+\llvert Y_s\rrvert
^2\bigr]\,ds \biggr\}^{p/2}
\\
&&\quad\qquad{}  + (1+\varepsilon
)^{{p/2}+1}c_1^p \biggl[\int_0^T
\llvert Z_s\rrvert ^2\,ds \biggr]^{p/2}
\\
&&\qquad \le C_{p,\varepsilon } \biggl\{ \biggl[\int_0^T
\bigl\llvert \sigma (s, 0)\bigr\rrvert ^2\,ds \biggr]^{p/2} +
T^{{p/2}-1}\int_0^T \bigl[\llvert
X_s\rrvert ^p+\llvert Y_s\rrvert
^p\bigr]\,ds \biggr\}
\\
&&\quad\qquad{}  + (1+\varepsilon )^{{p/2}+1}c_1^p
\biggl[\int_0^T \llvert Z_s
\rrvert ^2\,ds \biggr]^{p/2}.
\end{eqnarray*}
In the above, $\varphi (s,0)\stackrel{\triangle}{=}\varphi
(s,0,0,0)$, for $\varphi =b,\sigma $, respectively. Then
%
\begin{eqnarray}
\label{XLp} \mathbb{E}\bigl[\llvert X_t\rrvert ^p
\bigr] &\le& C_{p,\varepsilon } \bigl[\llvert x\rrvert ^p+I_p^p
+ \delta ^{p/2}\llVert \Theta \rrVert _{w,p}^p
\bigr]\nonumber
\\
&&{}+ (1+\varepsilon )\psi_1(p)^p \bigl[C_{p,\varepsilon }
\bigl[I_p^p + \delta ^{{p/2}-1}\llVert \Theta \rrVert
_{w,p}^p\bigr]\nonumber
\\
&&{} + (1+\varepsilon )^{{p/2}+1}c_1^p
\llVert \Theta \rrVert _{w,p}^p \bigr]
\\
&\le& C_{p,\varepsilon } \bigl[\llvert x\rrvert ^p+I_p^p
+ \delta ^{{p/2}-1}\llVert \Theta \rrVert _{w,p}^p
\bigr]\nonumber
\\
&&{} + (1+\varepsilon )^{{p/2}+2}\psi _1(p)^pc_1^p
\llVert \Theta \rrVert _{w,p}^p.\nonumber
\end{eqnarray}

Next, by Theorem \ref{teo-wellposed} we have
\begin{eqnarray*}
\llvert Y_t\rrvert ^2 \le C \mathbb{E}_t
\biggl[ \llvert X_t\rrvert ^2 + \bigl\llvert g(0)\bigr
\rrvert ^2 + \biggl(\int_t^T \bigl[
\llvert b\rrvert +\llvert f\rrvert \bigr](s,0)\,ds \biggr)^2 + \int
_t^T \bigl\llvert \sigma (s,0)\bigr\rrvert
^2\,dt \biggr].
\end{eqnarray*}
This implies that
%
\begin{eqnarray}
\label{YLp} \mathbb{E}\bigl[\llvert Y_t\rrvert ^p
\bigr] &\le& C_p\mathbb{E}\bigl[\llvert X_t\rrvert
^p\bigr] + C_pI_p^p.
\end{eqnarray}
In particular,
%
\begin{eqnarray}
\label{YLp0} \llvert Y_0\rrvert ^p &\le&
C_p\bigl[\llvert x\rrvert ^p + I_p^p
\bigr].
\end{eqnarray}

Moreover, following standard arguments
\begin{eqnarray*}
&& \mathbb{E} \biggl[\biggl\llvert \int_0^T
Z_t \,dB_t\biggr\rrvert ^p \biggr]
\\
&&\qquad = \mathbb{E}
\biggl[\biggl\llvert g(X_T) - g(0) + g(0) - Y_0 + \int
_0^T f(t, \Theta _t)\,dt\biggr\rrvert
^p \biggr]
\\
&&\qquad \le (1+\varepsilon ) \mathbb{E} \bigl[\bigl\llvert g(X_T) - g(0)
\bigr\rrvert ^p \bigr]
\\
&&\quad\qquad{} + C_{p,\varepsilon }\mathbb{E} \biggl[\bigl
\llvert g(0)\bigr\rrvert ^p + \llvert Y_0\rrvert
^p + \biggl\llvert \int_0^T f(t,
\Theta _t)\,dt\biggr\rrvert ^p \biggr]
\\
&&\qquad \le (1+\varepsilon )c_3^p \mathbb{E} \bigl[\llvert
X_T\rrvert ^p \bigr] + C_{p,\varepsilon } \bigl[ \llvert x
\rrvert ^p + I_p^p + \delta ^{p/2}
\llVert \Theta \rrVert _{w,p}^p \bigr].
\end{eqnarray*}
Now by the second inequality in (\ref{MLp}) and (\ref{XLp}), we have
%
\begin{eqnarray}\label{ZLp}
\qquad && \mathbb{E} \biggl[ \biggl(\int_0^T \llvert Z_t\rrvert ^2 \,dt \biggr)^{p/2} \biggr]\nonumber
\\
&&\qquad   \le    (1+\varepsilon
)c_3^p \bigl\llvert \psi_2(p)\bigr\rrvert
^p \mathbb{E} \bigl[\llvert X_T\rrvert ^p
\bigr] + C_{p,\varepsilon } \bigl[ \llvert x\rrvert ^p +
I_p^p + \delta ^{p/2}\llVert \Theta \rrVert
_{w,p}^p \bigr]
\nonumber\\[-8pt]\\[-8pt]\nonumber
&&\qquad  \le   (1+\varepsilon
)^{{p/2}+3}\bigl[\psi(p)c_1c_3\bigr]^p
\llVert \Theta \rrVert _{w,p}^p\nonumber
\\
&&\quad\qquad{}  + C_{p,\varepsilon } \bigl[
\llvert x\rrvert ^p+I_p^p + \delta
^{{p/2}-1}\llVert \Theta \rrVert _{w,p}^p \bigr].\nonumber
\end{eqnarray}
Set $\varepsilon =1$ in (\ref{XLp}), and plug (\ref{XLp}), (\ref
{YLp}), (\ref{ZLp}) into (\ref{ThLp0}), we get
%
\begin{eqnarray}
\label{ThLp} \llVert \Theta \rrVert _{w,p}^p
  &  \le  &  \mathbb{E} \biggl[
\biggl(\int_0^T \llvert Z_t\rrvert
^2 \,dt \biggr)^{p/2} \biggr] + \delta \sup
_{0\le t\le T} \mathbb{E}\bigl[\llvert X_t\rrvert
^p + \llvert Y_t\rrvert ^p\bigr]\nonumber
\\
  &  \le  &  \bigl[(1+
\varepsilon )^{{p/2}+3}\bigl[\psi(p)c_1c_3
\bigr]^p + C_{p,\varepsilon
}\delta ^{{p/2}-1}
\\
&&{} + C_p \delta \bigr]\llVert \Theta \rrVert _{w,p}^p +
C_{p,\varepsilon }\bigl[\llvert x\rrvert ^p+I_p^p
\bigr].\nonumber
\end{eqnarray}
Denote
\begin{eqnarray*}
c_p \stackrel{\triangle} {=}\bigl[\psi(p)c_1c_3
\bigr]^p <1.
\end{eqnarray*}
We may first choose $\varepsilon $ such that $(1+\varepsilon
)^{{p/2}+3}[\psi(p)c_1c_3]^p = {2c_p + 1\over3}$, and then
choose $\delta $ such that $C_{p,\varepsilon }\delta ^{{p/2}-1}
+ C_p\delta = {1-c_p\over6}$. Then (\ref{ThLp}) implies that
\begin{eqnarray*}
\llVert \Theta \rrVert _{w,p}^p \le{c_p+1\over2}
\llVert \Theta \rrVert _{w,p}^p + C_p\bigl[
\llvert x\rrvert ^p+I_p^p\bigr].
\end{eqnarray*}
Since ${c_p+1\over2}<1$, we obtain
$\llVert  \Theta \rrVert  _{w,p}^p \le C_p[\llvert  x\rrvert  ^p+I_p^p]$.
Now following standard arguments we can prove (\ref{Lp}) straightforwardly.

Finally, the claim on $(\hat Y, \hat Z)$ follows from Theorem \ref
{teo-linearuniq} immediately.
\end{pf}

We note that if $\sigma =\sigma (t,x,y)$, then we could simply take
$c_1=0$. Note that $\psi^{-1}(\infty) = \infty$, by combining the
arguments in Theorems \ref{teo-wellposed1} and \ref{teo-Lp} [noting~(\ref{p})], we obtain the following result immediately.

\begin{cor}
\label{cor-Lp} Let $(b,\sigma,f,g)$ satisfy the conditions in
Theorem \ref{teo-wellposed1}. For any $p\ge2$, if $I_p<\infty$, then
the unique solution $\Theta $ of FBSDE (\ref{FBSDE}) is in $L^p$ and
satisfies (\ref{Lp}).
Consequently, the corresponding
characteristic BSDE (\ref{hatYBSDE}) has a unique solution
$(\hat Y, \hat Z)$
satisfying (\ref{hatYproperty}) and (\ref{hatZproperty}).
\end{cor}

For FBSDE (\ref{linearFBSDE}), we have $I_p=0$ for all $p\ge2$, which
leads to the following result.

\begin{cor}
\label{cor-Lplinear} Assume the linear FBSDE (\ref{linearFBSDE})
satisfy the conditions in
Theorem \ref{teo-wellposed} (or Theorem \ref{teo-wellposed1}). Then
any $2\le p< \psi^{-1}({1\over c_1c_3})$, the unique solution~$\Theta
$ of FBSDE (\ref{linearFBSDE}) is in $L^p$. Consequently, the corresponding
characteristic BSDE (\ref{hatYBSDE}) has a unique solution
$(\hat Y, \hat Z)$
satisfying (\ref{hatYproperty}) and (\ref{hatZproperty}).
\end{cor}

Finally, as an application of Corollary \ref{cor-Lplinear}, we prove
the comparison theorem.

\begin{teo}[(Comparison)]
\label{teo-comparison} Assume both $(b, \sigma, f, g)$ and $(b,
\sigma,
\tilde f, \tilde g)$ satisfy the same conditions (i.e., they belong
to the same case) in Theorem
\ref{teo-wellposed} (or Theorem \ref{teo-wellposed1}), and let $u,
\tilde u$ be the
corresponding random fields. If $f\le\tilde f, g\le\tilde g$, then
$u\le\tilde u$.
\end{teo}

\begin{pf}
Without loss of generality, we shall prove the result only
at $t=0$. Let $\Theta, \tilde\Theta \in L^2$ be the
corresponding solutions to the FBSDE (\ref{FBSDE}) associated to
$(b,\sigma,f,g)$ and $(b,\sigma,\tilde f,\tilde g)$, respectively.
Denote
$\Delta \Theta _t \stackrel{\triangle}{=}\Theta _t-\tilde\Theta
_t$, and define $\varphi _i$ similar to (\ref{tildephi}) for
$\varphi =b$, $\sigma $, $f$, respectively. Then $\Delta \Theta $
would be the unique solution
to the following linear FBSDE:
%
\begin{eqnarray}
\label{FBSDE-comparison} \cases{ \displaystyle\Delta X_t= \int
_0^t (b_1\Delta
X_s+b_2\Delta Y_s+b_3\Delta Z_s)\,ds
\vspace*{5pt}\cr
\displaystyle \phantom{\Delta X_t=}{} +\int_0^t (\sigma
_1\Delta X_s+\sigma_2\Delta
Y_s+\sigma _3\Delta Z_s)\,dB_s;
\vspace*{5pt}\cr
\displaystyle \Delta Y_t= h\Delta X_T+\Delta g(\tilde X_T)
\vspace*{5pt}\cr
\displaystyle \phantom{\Delta Y_t=}{} + \int_t^T
\bigl(f_1\Delta X_s+f_2\Delta
Y_s+f_3\Delta Z_s+\Delta f(t,\tilde\Theta
_t)\bigr)\,ds
\vspace*{5pt}\cr
\displaystyle \phantom{\Delta Y_t=}{}  -\int_t^T\Delta
Z_s\,dB_s.}
\end{eqnarray}

Let $(\hat Y, \hat Z)$ denote the unique solution to BSDE (\ref
{hatYBSDE}) which, by Corollary \ref{cor-Lplinear}, satisfies (\ref
{hatYproperty}) and (\ref{hatZproperty}).
Denote
\begin{eqnarray*}
\delta Y &\stackrel{\triangle} {=}&\Delta Y- \hat Y \Delta X,
\\
\delta Z &\stackrel{\triangle} {=}&\Delta Z - \hat Z \Delta X - \hat Y[\sigma _1
\Delta X+\sigma _2\Delta Y+\sigma _3\Delta Z],
\end{eqnarray*}
and define $\beta, \gamma $ and $\Gamma $ by (\ref{betagamma}) and
(\ref{Gamma}). Applying It\^{o}'s formula, we have
\begin{eqnarray*}
\delta Y_0 = \Gamma _0 \delta Y_0 =
\Gamma _T \Delta g(\tilde X_T) + \int_0^T
\Gamma _t \Delta f(t, \tilde\Theta _t)\,dt - \int
_0^T \Gamma _t[\gamma
_t \delta Y_t + \delta Z_t]\,dB_t.
\end{eqnarray*}
Now by (\ref{hatZproperty}) and following similar arguments as in
Theorem \ref{teo-linearuniq} one can easily show that $ \int_0^t
\Gamma _s[\gamma _s \delta Y_s + \delta Z_s]\,dB_s$ is a true
martingale. Then by our assumptions we see that
\begin{eqnarray*}
u(0,x) - \tilde u(0,x) = \Delta Y_0 = \delta Y_0 =
\mathbb{E} \biggl\{ \Gamma _T \Delta g(\tilde X_T) + \int
_0^T \Gamma _t \Delta f(t, \tilde
\Theta _t)\,dt \biggr\}\le0.
\end{eqnarray*}
This proves the theorem.
\end{pf}

\begin{rem}
We notice that we cannot get $\Delta Y_t\ge0$ even $\Gamma
_t\ge0$,
$0\le t\le T$, in the above proof. This coincides with the
results in Wu and Xu \cite{wx} (Theorem~3.2 and Counterexample~3.1). However, for the corresponding random decoupling field, the
comparison theorem holds over all time which coincides with Theorem
4.1 in Cvitanic and Ma \cite{CM} by virtue of PDE method under
Markovian frame work.
\end{rem}

\begin{appendix}\label{app}
\section*{Appendix}
\label{sect-appendix}
\setcounter{equation}{0}

In this Appendix, we complete the technical proofs for some results in
Section~\ref{sec5}.

\begin{pf*}{Proof of Lemma \ref{lem-ODEcomp}}
We first show the
existence. Define a truncation function
\begin{eqnarray*}
\tilde F(t, y) \stackrel{\triangle} {=}F \bigl(t, \by^1_t
\vee y \wedge\by^2_t \bigr),
\end{eqnarray*}
then by assumption (iii) $\tilde F$ is uniformly Lipschitz
continuous in $y$ with a Lipschitz constant $L$, and thus the following
ODE has a unique solution $\tilde\by$:
%
\begin{eqnarray}
\label{ODEcomp-tildey} \tilde\by_t = h + \int_t^T
\tilde F(s,\tilde\by_s)\,ds, \qquad t\in[0,T].
\end{eqnarray}
We claim that
%
\begin{eqnarray}
\label{y1<y<y2} \by^1 \le\tilde\by\le\by^2.
\end{eqnarray}
This would lead to that $\tilde F(t, \tilde\by_t) = F(t,\tilde\by
_t)$. Thus, $\tilde\by$ is a solution to ODE (\ref{ODE0}) and~(\ref
{y1<y<y2}) holds.

In fact, denote $\Delta \by^2 \stackrel{\triangle}{=}\by^2 -
\tilde\by$, $\Delta h^2\stackrel{\triangle}{=}h^2 - h$, $\Delta
F^2 \stackrel{\triangle}{=}F^2 - F$. Note that
$ F(t, \by^2_t) = \tilde F(t, \by^2_t)$, we have
\begin{eqnarray*}
\Delta \by^2_t &=& \Delta h^2 +
C^2 + \int_t^T
\bigl[F^2\bigl(s, \by^2_s\bigr) - \tilde F(s,
\tilde\by_s) - c^2_s\bigr]\,ds
\\
&=& \Delta h^2 + C^2 + \int_t^T
\bigl[\Delta F^2\bigl(s, \by^2_s\bigr) +
\alpha _s \Delta \by^2_s -
c^2_s\bigr]\,ds,
\end{eqnarray*}
where
$ \alpha _s \stackrel{\triangle}{=}{\tilde F(s, \by^2_s)-\tilde
F(s, \tilde\by_s)\over\Delta \by^2_s}\textbf{1}_{\{\Delta \by
^2_s\neq0\}}$ satisfies $\llvert  \alpha \rrvert  \le L$. Now define $\gamma _t
\stackrel{\triangle}{=}\break \exp(\int_0^t \alpha _s \,ds)>0$. Then
\begin{eqnarray*}
\gamma _t \Delta \by^2_t &=& \gamma
_T\bigl[\Delta h^2 + C^2\bigr] + \int
_t^T \gamma _s\bigl[\Delta
F^2\bigl(s, \by^2_s\bigr) -
c^2_s\bigr]\,ds
\\
&=& \gamma _T\Delta h^2 + \int_t^T
\gamma _s\Delta F^2\bigl(s, \by ^2_s
\bigr)\,ds + \gamma _T \biggl[C^2 - \int
_t^T \gamma _T^{-1}\gamma
_s c^2_s\,ds \biggr]\ge0.
\end{eqnarray*}
This implies that $\tilde\by\le\by^2$. Similarly, we have $\tilde
\by\ge\by^1$.

It remains to prove the uniqueness. Let $\by$ be an arbitrary solution
to ODE (\ref{ODE0}) satisfying (\ref{y1<y<y2}). Then $\tilde F(t, \by
_t) = F(t,\by_t)$, and thus $\by$ satisfies ODE (\ref
{ODEcomp-tildey}). By the uniqueness of ODE (\ref{ODEcomp-tildey}) we
have $\by= \tilde\by$, and thus uniqueness follows.
\end{pf*}

\begin{pf*}{Proof of Theorem \ref{teo-constant-F1}}
(\textit{Necessity}). For simplicity, let us rewrite (\ref{F1}) as
%
\begin{eqnarray}
\label{a123} F(y)=f_1+a_1 y+a_2y^2+a_3y^3,
\end{eqnarray}
where $a_3=\sigma _2b_3$, $a_2=b_2+f_3\sigma _2+b_3\sigma _1$,
$a_1=f_2+b_1+\sigma _1f_3$.

We shall show that if none of (i)--(iii) holds, then the solution of
ODE (\ref{ODE1}) will blow-up in finite time, which would complete
the proof. To this end, we assume without loss of generality that
$F(h)\ge0$. [The case when $F(h)\le0$ can be argued in the same way
but using the conditions (ii) and (iii).] Since (i) does not hold, $F$ has
no zero point in $[h, \infty)$, and hence $F(h)>0$. Now since (iii)
does not hold, $\llvert  a_3\rrvert  +\llvert  a_2\rrvert  \neq0$.
Note that if $a_3 < 0$ or $a_3 =
0$ but $a_2 <0$, then $\lim_{y\to\infty}F(y) = -\infty$ which,
together with $F(h)>0$, will imply that $F$ has a zero point in $[h,
\infty)$, a contradiction. Thus, we need only check the case where
\textit{either} ``$a_3>0$'' \textit{or} ``\mbox{$a_3=0$}, \mbox{$a_2>0$}.''
We investigate the two
cases separately.

\begin{longlist}
\item[\textit{Case} 1.] Assume $a_3>0$. We claim that there exist $\varepsilon
>0$ and $y_1< h$ such that
%
\begin{eqnarray}
\label{a3} F(y)\ge\varepsilon (y-y_1)^{3}\qquad\mbox{for all } y\ge h.
\end{eqnarray}
Indeed, in this case $F(y)$ is a polynomial of degree 3, it must
have at least one real zero point. By our assumption, $F$ has no zero
point after $h$, then all real zero points must be in $(-\infty,
h)$. If there are three real zero points (possibly equal), we list them as
$-\infty<y_1\le y_2\le y_3<h$. Then for any $y\ge h$, one has
%
\begin{eqnarray}
\label{a31} F(y)=a_3\prod_{i=1}^3
(y-y_i)\ge a_3(y-y_1)^3.
\end{eqnarray}
On
the other hand, if $F$ has only one real zero point, denoted as
$y_1$, then we may write
\begin{eqnarray*}
F(y)=a_3(y-y_1) \bigl[(y-y_2)^{2}+c
\bigr]\qquad\mbox{for some }c>0.
\end{eqnarray*}
Note that the function $\tilde F(y) \stackrel{\triangle}{=}
a_3[(y-y_2)^{2}+c] (y-y_1)^{-2}$ is continuous for $y\in[h, \infty)$,
$\tilde F(y) >0$ and $\lim_{y\to\infty} \tilde F(y) = a_3>0$. Then
\[
\varepsilon \stackrel{\triangle} {=}\inf_{y\ge h}
{a_3[(y-y_2)^{2}+c]\over(y-y_1)^2}>0.
\]
Thus, noting that $y-y_1>0$ for $y\ge h$,
\[
F(y)=a_3(y-y_1) \bigl[(y-y_2)^{2}+c
\bigr]\ge\varepsilon (y-y_1)^3\qquad\mbox{for all } y\ge h.
\]
This, together with (\ref{a31}), proves (\ref{a3}).

Now consider the following ODE:
%
\begin{equation}
\label{ODE11} \tilde\by_t=h+\int_t^T
\varepsilon (\tilde\by_t-y_1)^{3}\,dt.
\end{equation}
%
Solving this ODE, we have
$\tilde\by_t-y_1={1\over\sqrt{2\varepsilon (t-T)+(h-y_1)^{-2}}}$.
Thus, if $T> {1\over2\varepsilon (h-y_1)^2}$, then the solution
$\tilde\by_t$ blows up at $t=T- {1\over2\varepsilon (h-y_1)^2}\in
(0,T)$. On the other hand,
by comparison theorem we can easily show that $\by_t\ge\tilde\by_t$.
Thus, the solution of (\ref{ODE1}) will blow-up at finite time as well.

\item[\textit{Case} 2.] Assume $a_3=0$ and $a_2>0$. Following similar arguments,
in this case we have
$ F(y)\ge\varepsilon (y-y_1)^{2}$, for all $ y\ge h$,
and similarly $\by$ will blow up if $T$ is large enough.\quad\qed
\end{longlist}\noqed
\end{pf*}

\begin{pf*}{Proof of Theorem \ref{teo-constant-F3}}
(\textit{Necessity}):
\begin{longlist}[(iii)]
\item[(i)] Assume $h< \sigma _3^{-1}$, $F(h)\le0$, and $\alpha _3 \stackrel
{\triangle}{=}b_2 - b_3\sigma _2 \sigma _3^{-1} \neq0$. We show
that either $F$ has a zero in $(-\infty, h]$ or $\by$ blows up when
$T$ is large enough.

Indeed, if $\alpha _3>0$, then $\lim_{y\to-\infty} F(y) =\infty$. Note
that $F$ is continuous for $y\in(-\infty, h]$. These, together with
$F(h)\le0$, imply that $F$ has a zero point in $(-\infty, h]$. We
now assume $\alpha _3 <0$. Denote $\tilde F(y) \stackrel{\triangle
}{=}-{F(y)\over
(h+1-y)^2}$. In $(-\infty, h]$, if $F$ has no zero point, then
$\tilde F$ is continuous, has no zero point, and $\lim_{y\to-\infty}
\tilde F(y) = -\alpha _3>0$. Denote $\varepsilon \stackrel{\triangle
}{=}\inf_{y\le h} \tilde
F(y)>0$. Then we have
\[
F(y) \le-\varepsilon (h+1-y)^2\qquad\mbox{for all } y\le h.
\]
Following the arguments for the proof of the necessary part of Theorem
\ref{teo-constant-F1}, we prove that $\by$ blows up when $T$ is large.

\item[(ii)] Assume $h> \sigma _3^{-1}$, $F(h)\ge0$, and $\alpha _3 \neq0$.
Similarly, we can show that either $F$ has a zero point in $[h,
\infty)$ or $\by$ blows up when $T$ is large enough.

\item[(iii)] Assume $h< \sigma _3^{-1}$ and $F(h)\ge0$. We show that either
$F$ has a zero point in $[h,\sigma _3^{-1})$ or $\by$ violates
(\ref{yproperty}) when $T$ is large enough.

Indeed, recall the $\alpha _0$ in (\ref{F2}). If $\alpha _0 <0$, then
$\lim_{y\uparrow\sigma _3^{-1}} F(y) = -\infty$. This implies that $F$
has a zero point in $[h,\sigma _3^{-1})$.

If $\alpha _0>0$ and $F$ has no zero point in $[h,\sigma _3^{-1})$.
Denote $\tilde F(y) \stackrel{\triangle}{=}F(y) [\sigma _3^{-1}-y]$.
Then in $[h, \sigma _3^{-1})$, $\tilde F$ is continuous, $\tilde F>0$,
and $\lim_{y\uparrow\sigma _3^{-1}} \tilde F(y) = \alpha _0>0$. Denote
\begin{eqnarray*}
\varepsilon \stackrel{\triangle} {=}\inf_{y\in[h, \sigma _3^{-1})} \tilde F(y)
>0\quad\mbox{and thus}\quad F(y) \ge\varepsilon \bigl(\sigma _3^{-1}-y
\bigr)^{-1}\qquad\mbox{for }y\in\bigl[h, \sigma _3^{-1}\bigr).
\end{eqnarray*}
Let $\tilde\by$ solve the following ODE:
\begin{eqnarray*}
\tilde\by_t = h + \int_t^T
\varepsilon \bigl(\sigma _3^{-1}-\tilde\by _s
\bigr)^{-1} \,ds,
\end{eqnarray*}
we obtain explicitly
$(\sigma _3^{-1} -\tilde\by_t)^2 = (\sigma _3^{-1} - h)^2 -
2\varepsilon (T-t)$.
Let $T \ge{1\over2\varepsilon } (\sigma _3^{-1} - h)^2$. Then for
$t = T - {1\over2\varepsilon } (\sigma _3^{-1} - h)^2\in[0,T]$, we
have\vspace*{1pt} $\tilde\by_t = \sigma _3^{-1}$. By comparison, we see that
$(1-\sigma _3\by)^{-1}$ would blow up.

Finally, if $\alpha _0=0$ and $F$ has no zero point in $[h,\sigma _3^{-1})$.
Then $F$ is continuous and positive on $[h, \sigma _3^{-1}]$. Denote
$\varepsilon
\stackrel{\triangle}{=}\inf_{y\in[h, \sigma _3^{-1}]} F(y) >0$,
and define
$\tilde
\by_t \stackrel{\triangle}{=}h + \int_t^T \varepsilon \,ds=h +
\varepsilon (T-t)$, $t\in[0,T]$.
Thus, if $T \ge\varepsilon ^{-1}[\sigma _3^{-1}-h]$, then $\tilde
\by_t =
\sigma _3^{-1}$ at $t=T - \varepsilon ^{-1}[\sigma _3^{-1}-h]$. By
comparison again,
we see that $(1-\sigma _3\by)^{-1}$ would blow up.

\item[(iv)] Assume $h>\sigma _3^{-1}$ and $F(h)\le0$. We can similarly show
that either $F$ has a zero point in $(\sigma _3^{-1}, h]$ or $\by$
violates (\ref{yproperty}) when $T$ is large enough.\quad\qed
\end{longlist}\noqed
\end{pf*}
\end{appendix}

\section*{Acknowledgments}
Part of this work was completed while Zhen Wu and
Detao Zhang were visiting the Department of
Mathematics, University of Southern California, whose hospitality is
greatly appreciated.



\printaddresses
\end{document}